 
\input amstex
\documentstyle{amsppt}

\font \aufont=  cmr10 
\font\titfont=  cmr12

\magnification=1200

\NoBlackBoxes

\pagewidth{5.8in}
\pageheight{7.8in}






 \define\zz{ {{\bold{Z}_2}}}


 \define\calc{\Cal C}
 
 \define\calj{\Cal J}
 \define\calk{\Cal K}
 
 \define\calq{\Cal Q}
 \define\calr{\Cal R}
 \define\calz{\Cal Z}



 \define\cycd#1#2{{\calz}_{#1}(#2)}
 \define\cyc#1#2{{\calz}^{#1}(#2)}
 \define\cych#1{{{\calz}^{#1}}}
 \define\cycp#1#2{{\calz}^{#1}(\bbp(#2))}
 \define\cyf#1#2{{\cyc{#1}{#2}}^{fix}}
 \define\cyfd#1#2{{\cycd{#1}{#2}}^{fix}}

 \define\crl#1#2{{\calz}_{\bbr}^{#1}{(#2)}}

 \define\crd#1#2{\widetilde{\calz}_{\bbr}^{#1}{(#2)}}

 \define\crld#1#2{{\calz}_{#1,\bbr}{(#2)}}
 \define\crdd#1#2{\widetilde{\calz}_{#1,{\bbr}}{(#2)}}

 \define\crlh#1{{\calz}_{\bbr}^{#1}}
 \define\crdh#1{{\widetilde{\calz}_{\bbr}^{#1}}}
 \define\cyav#1#2{{{\cyc{#1}{#2}}^{av}}}
 \define\cyavd#1#2{{\cycd{#1}{#2}}^{av}}

 \define\cyaa#1#2{{\cyc{#1}{#2}}^{-}}
 \define\cyaad#1#2{{\cycd{#1}{#2}}^{-}}

 \define\cyq#1#2{{\calq}^{#1}(#2)}
 \define\cyqd#1#2{{\calq}_{#1}(#2)}

 \define\cqt#1#2{{\calz}_{\bbh}^{#1}{(#2)}}
 \define\cqtav#1#2{{\calz}^{#1}{(#2)}^{av}}
 \define\cqtrd#1#2{\widetilde{\calz}_{\bbh}^{#1}{(#2)}}

 \define\cyct#1#2{{\calz}^{#1}(#2)_\zz}
 \define\cyft#1#2{{\cyc{#1}{#2}}^{fix}_\zz}
\define\cxg#1#2{G^{#1}_{\bbc}(#2)}
 \define\reg#1#2{G^{#1}_{\bbr}(#2)}

 \define\cyaat#1#2{{\cyc{#1}{#2}}^{-}_\zz}

 \define\pc{\bbp_{\bbc}}
 \define\pr{\bbp_{\bbr}}


 \define\fflag#1#2{{#1}={#1}_{#2} \supset {#1}_{{#2}-1} \supset
 \ldots \supset {#1}_{0} }
 \define\csus{\Sigma \!\!\!\! /\,}
 \define\vect#1{ {\Cal{V}ect}_{#1}}
 \define\BU{{BU}}
 \define\BO{{BO}}


 \define\chv#1#2#3{{\calc}^{#1}_{#2}(#3)}
 \define\chvd#1#2#3{{\calc}_{#1,#2}(#3)}
 \define\chm#1#2{{\calc}_{#1}(#2)}



 \define\Claim#1{\subheading{Claim #1}}

 \define\xrightarrow#1{\overset{#1}\to{\rightarrow}}

 \define\cx{\gamma}
 \define\cxx{\Gamma}

 \define\Real{\phi}
 \define\Reall{\Phi}

 \define\pont{j}
 \define\whit{\widetilde{j}}


\font \fr = eufm10

\hfuzz1pc 


\define\gM{\text{\fr M}}

\define\gR{\text{\fr R}}
\define\gK{\text{\fr K}}
\define\gu{\text{\fr u}}

\define\go{\text{\fr o}}

\define\gZ{\text{\fr Z}}

\define\bbz{\Bbb Z}

\define\bbr{\Bbb R}
\define\bbc{\Bbb C}
\define\bbh{\Bbb H}
\define\bbp{\Bbb P}

\define\bp{\bold P}
\define\bz{\bold Z}

\define\bb{\bold b}

\define\br{\bold r}

\define\cu{\Cal U}

\define\cc{\Cal C}

\define\cl{\Cal L}
\define\cg{\Cal G}

\define\cf{\Cal F}

\define\cz{\Cal Z}
\define\co#1{\Cal O_{#1}}
\define\ct{\Cal T}
\define\ci{\Cal I}


\define\a{\alpha}
\redefine\b{\beta}

\define\r{\rho}

\define\z{\zeta}
\define\x{\xi}
\define\la{\lambda}

\redefine\D{\Delta}
\define\G{\Gamma}

\define\p#1{{\bbp}^{#1}}

\define\equdef{\overset\text{def}\to=}

\define\blbx{\hfill  $\square$}
\redefine\qed{\blbx}

\define\pf{\subheading{Proof}}
\define\Lemma#1{\subheading{Lemma #1}}
\define\Theorem#1{\subheading{Theorem #1}}
\define\Prop#1{\subheading{Proposition #1}}
\define\Cor#1{\subheading{Corollary #1}}
\define\Note#1{\subheading{Note #1}}
\define\Def#1{\subheading{Definition #1}}
\define\Remark#1{\subheading{Remark #1}}
\define\Ex#1{\subheading{Example #1}}
\define\arr{\longrightarrow}

\define\Hom{\text{Hom}}

\redefine\Xi{X_{\infty}}

\define\jac#1#2{\left(\!\!\!\left(
\frac{\partial #1}{\partial #2}
\right)\!\!\!\right)}
\define\restrict#1{\left. #1 \right|_{t_{p+1} = \dots = t_n = 0}}

\define\SP#1#2{{\roman SP}^{#1}(#2)}

\define\coc#1#2#3{\cc^{#1}(#2;\, #3)}
\define\zoc#1#2#3{\cz^{#1}(#2;\, #3)}
\define\zyc#1#2#3{\cz^{#1}(#2 \times #3)}

\define\ar#1{\overset{#1}\to{\longrightarrow}}

\define\th#1#2{{\Bbb H}^{#1}(#2)}
\define\hth#1#2{\widehat{\Bbb H}^{#1}(#2)}


\define\tcz{\widetilde{\cz}}

\define\tf{\widetilde f}

\define\bad#1#2{\cf_{#2}(#1)}

\define\pch#1{\bbp_{\bbc}(\bbh^{#1})}

\def\l{\ell}

\def\I#1#2{I_{#1, #2}}
\def\Iav#1#2{I_{#1, #2}^{av}}
\def\GR#1#2{I_{#1, #2}}
\def\GRav#1#2{I_{#1, #2}^{av}}
\def\GRtil#1#2{\widetilde{I}_{#1, #2}}

\def\A{A}
\def\BLLMM{BLLMM}
\def\CW{CW}
\def\DT{DT}

\def\FM{FM}
\def\FL{FL}
\def\Fu{Fu}
\def\Lam{Lam}
\def\L{L}
\def\Li{Li}
\def\LLM{LLM}
\def\LM{LM}
\def\LMS{LMS}
\def\May{M}
\def\MS{MS}
\def\Mo{Mo}
\def\Se{Seg}

\def\tf{\widetilde{f}}
\def\z2t{\text{$\bbz_2\ct$}}

\def\dS{dS}
\def\D{D}

\def\bbk{{\Bbb K}}

\define\tc{\widetilde{c}}

\ \bigskip
\centerline{\bf \titfont  ALGEBRAIC CYCLES AND THE CLASSICAL 
GROUPS }

\vskip .12in

\centerline{\bf \titfont  Part I,\ \ \ Real Cycles }

\vskip .12in
\centerline{by}
\vskip .2in

\centerline{\bf  \aufont  H. Blaine Lawson, Jr., \ 
Paulo Lima-Filho,\ and\ 
Marie-Louise Michelsohn}

\vskip .3in
 \font\abfont=cmr10
 \topmatter
\abstract {\abfont
The groups of algebraic cycles on complex projective space 
$\bbp(V)$ are
known to have  beautiful and surprising properties.
Therefore, when  $V$ carries a real structure,  it is natural to ask for
the  properties of the groups of real algebraic cycles on $\bbp(V)$. 
Similarly, if  $V$ carries a quaternionic structure, one can  define 
quaternionic algebraic cycles and  ask
the same question. In this paper and its sequel 
  the homotopy structure of these cycle groups is completely 
determined. 
It  turns out to be quite  simple and to bear a direct
relationship to characteristic classes for the  classical groups.

It is   shown, moreover, that certain functors in $K$-theory extend
directly to  these groups.   
It is also shown that, after taking
colimits over dimension and codimension, the groups of real and
quaternionic  cycles   carry $E_{\infty}$-ring structures, and that the
maps extending the  $K$-theory functors are  $E_{\infty}$-ring maps. 
This gives a wide  generalization of the results in [\BLLMM] on the
Segal question.

The ring structure on the homotopy groups of these stabilized spaces
is explicitly computed.  In the real case it is a simple quotient
of a polynomial algebra on two generators corresponding to the first
Pontrjagin and first Stiefel-Whitney classes.

These calculations yield  an interesting
total characteristic class for real bundles.  It is a
mixture of integral and mod 2 classes  and   has   nice
multiplicative properties. The class is shown to be related to the
$\bbz_2$-equivariant Chern class on Atiyah's $KR$-theory.} 
\endabstract
\endtopmatter

\vskip .1in       
        
\centerline{Table of Contents}
\medskip  
       
\hskip 1in  1.  Introduction.
 
\hskip 1in  2.  Spaces of complex cycles.

\hskip 1in  3.  Spaces of real cycles.

\hskip 1in  4.  The ring structure.

\hskip 1in  5.  Extending functors from $K$-theory.  

\hskip 1in  6.  Equivariant infinite loop space structures and $KR$-theory.

\hskip 1in  7.  A new total characteristic class.

\hskip 1in  8.  Proof of Theorems 3.3 and 3.4.

\hskip 1in  9.  Proof of Theorem 4.1.       

\hskip 1in  A.  Appendix:  Canonical splittings. 

\vfill\eject

\document

\subheading{\S 1.  Introduction}

In recent years a number of results have been proved about
the topological groups of algebraic cycles on an algebraic variety $X$
over $\bbc$. It has been shown for example that  when $X$ is
projective space, these groups provide useful models for  basic 
classifying spaces in algebraic topology and for certain  universal
characteristic maps between them.  They also yield 
certain new infinite loop space structures on products of
Eilenberg-MacLane spaces which make the total Chern class an 
infinite
loop map.  (See [\L$_2$] for a survey.)

Now  when   $X$ has a real
structure, it is natural to consider  the {\sl real algebraic
cycles} on  $X$.  These are simply the cycles defined over
$\bbr$, or equivalently, the cycles fixed by the Galois group 
$Gal(\bbc/\bbr)$. When $X$ is projective space $P(V)$ the set of real
cycles of codimension-$q$ forms a topological group
$\cz_{\bbr}^q(\bbp(V))$ whose homotopy-type is independent of 
$V$
[\Lam].  The first  main result of this paper  is the determination of
the topological structure of $\cz_{\bbr}^q(\bbp(V))$.  We show that it
canonically decomposes into a product of Eilenberg-MacLane spaces 
for
the groups $\bbz$ and $\bbz_2$. (See Theorem 3.3 below.)\ \  The
resulting structure is rather complicated when compared to the 
complex
case.

Our first explanation for the richness of this structure comes from
considering the colimit $\cz^{\infty}_{\bbr}$ of these groups over
dimension and codimension.  Here the algebraic join of cycles
induces a ring structure on the homotopy groups and we show that
as a ring
$$
\pi_* \cz^{\infty}_{\bbr} \ \cong\ \bbz[x,y]/(2y)
\tag{1.1}
$$
where $x$ corresponds to the generator of 
$\pi_4\cz^{\infty}_{\bbr}\cong \bbz$ and  $y$ 
corresponds to the generator of 
$\pi_1\cz^{\infty}_{\bbr}\cong \bbz_2$.

  Now the Grassmannian $\cg^q(\bbp(V))$ of codimension-$q$ planes
in $V$ includes naturally into  $\cz^q(\bbp(V))$ as degree-1
cycles.  Restricting to real points gives an inclusion
$\cg^q_{\bbr}(\bbp(V))  \to \cz_{\bbr}^q(\bbp(V))$ which stabilizes
to a mapping
$$
P : BO_q \ \arr\ \cz_{\bbr}^q(\bbp^{\infty}).
$$
This map represents an interesting total characteristic class
which, via  Theorem 3.3, is an explicit 
combination of integral and mod 2 cohomology classes and which has 
the
property that for real vector bundles $E$ and $F$ 
$$
P(E\oplus F) \ = \ P(E)P(F).
$$
 
In \S 6 we show that $ \cz^{\infty}_{\bbr}$ carries the structure of 
an
$E_{\infty}$-ring space and thus gives rise to an $E_{\infty}$-ring
spectrum. The additive deloopings in this spectrum are the standard
deloopings of Eilenberg-MacLane spaces.  The multiplicative 
deloopings 
extend the product in the group of multiplicative units of the theory.
We show that for the multiplicative deloopings, the limiting map 
$$
P: {\bold{BO}} \ \arr\  \cz^{\infty}_{\bbr}
$$ 
is an
infinite loop map yielding a map of spectra
$
P:\gK\go \to \gM_{\bbr}
$
from connective K-theory to the multiplicative spectrum of the 
theory.

The cycle groups admit two natural homomorphisms:
$$
\cz_{\bbc}^q(\bbp(V))
\ \longleftarrow\ 
\cz_{\bbr}^q(\bbp(V)) \ \arr\ \tcz_{\bbr}^q(\bbp(V)). 
$$
The left mapping  is the obvious inclusion.  The right mapping is
projection to the Galois quotient  $\tcz_{\bbr}^q(\bbp(V)) =
\cz_{\bbr}^q(\bbp(V))/\cz_{\text{av}}^q(\bbp(V))$ where
$\cz_{\text{av}}^q(\bbp(V)) =  \{c+\overline c : c \in
\cz_{\bbc}^q(\bbp(V))\}$. The colimits  of these spaces
have  $E_{\infty}$-ring structures for which the limiting maps
$$
 \cz_{\bbc}^{\infty} 
\ \longleftarrow\ 
\cz_{\bbr}^{\infty}  \ \arr\ \tcz_{\bbr}^{\infty}  
$$
are $E_{\infty}$-ring maps.  It is known that 
$\pi_*\cz_{\bbc}^{\infty} \cong \bbz[s]$ where $s$ corresponds to 
the
generator of $\pi_2\cz_{\bbc}^{\infty} \cong \bbz$ and 
$\pi_*\cz_{\bbr}^{\infty} \cong
\bbz_2[y]$ where $y$ corresponds to the generator of 
$\pi_1\cz_{\bbr}^{\infty} \cong \bbz_2$. Under the isomorphism 
(1.1) 
we show that the maps above induce ring homomorphisms
$$
 \bbz[s]
\ \longleftarrow\ 
\bbz[x,y]/(2y) \ \arr\ \bbz_2[y]
$$
given by $x\mapsto s^2$ and $y\mapsto y$.

Composing with the mapping $P$ gives two new mappings
$$\CD
\ @. \  @. \ @. \ @.  
\cz_{\bbc}^{\infty} \cong \prod_{k\geq0}K(\bbz,2k)   \\ 
\ @.\   @.  \nearrow @. \ @. \       \\ 
\text{\bf BO} @>{P}>> {\cz}^{\infty}_{\bbr} @. \ @. \ @. \  \\ 
\ @. \ @. \searrow @. \ @. \  \\
\ @. \  @. \ @. \ @.  \widetilde{\cz}^{\infty}_{\bbr} 
\cong\prod_{k\geq0}K(\bbz_2,k). \\
\endCD$$
The top composition classifies the total Chern class of the
complexification, and the bottom classifes the total Stiefel-Whitney
class. Thus the characteristic class $P$ carries all this information.
Furthermore, the maps above all extend to infinite loop maps.

Surprisingly other natural functors in K-theory, such as the forgetful
functor, extend from Grassmannians to the spaces of all cycles 
yielding
new proofs of relations between characteristic classes. (See \S 5.)
In \S 6 these maps  are also shown to be infinite loop maps.

There is a  unifying perspective on the results discussed above. 
For this we revisit the map
$$
c: {\bold{BU}} \ \arr\ \cz_{\bbc}^{\infty}
\tag{1.2}
$$
and recall that it is a $\bbz_2$-map with respect to complex
conjugation.  Thus we plunge into the world of $\bbz_2$-spaces,
$\bbz_2$-maps, and $\bbz_2$-equivariant homotopy theory.  Note 
that a
$\bbz_2$-space is just a Real space in the sense of Atiyah  [\A]. 
Furthermore,  ${\bold{BU}}$ is the classifying space for Atiyah's
$KR$-theory. We prove in \S 6 that $\cz_{\bbc}^{\infty}$ has the 
structure of
a   $\bbz_2 E_{\infty}$-ring space and that $c$ is a
$\bbz_2$-equivariant  infinite loop map into the multiplicative
structure.

In his thesis, Pedro dos Santos has proved that there is a canonical
$\bbz_2$-equivariant homotopy equivalence
$$
\cz_{\bbc}^{\infty}\ \cong\ \prod_{k\geq0} K(\underline{\bbz}, 
\bbr^{n,n})
\tag{1.3}
$$
where $K(\underline{\bbz}, \bbr^{n,n})$ denotes the equivariant 
Eilenberg-MacLane space classifying $\bbz_2$-equivariant 
cohomology
indexed at the representation  $\bbr^{n,n} = \bbc^n$ (with action 
given
by complex conjugation) and with coefficients in the constant Mackey
functor $\underline{\bbz}$.  
He furthermore shows that with respect to (1.3) the algebraic join
pairing classifies the equivariant cup product and the mapping 
(1.2) classifies the equivariant total Chern class in $KR$-theory.
 Our characteristic mapping $P$ represents
the restriction of this equivariant Chern class to the fixed-point sets.
(See \S 6 for details.) 

Analogous results for the quaternionic case are proved in Part II of
this paper.

The authors would like thank Pedro Santos and Daniel Dugger for 
several very useful remarks and conversations relating to this work. 
They would also like to mention that in his thesis, Jacob Mostovoy
computed some of the homotopy groups that appear in this paper. His
results were announced in [\Mo$_1$] and subsequently appeared in [\Mo$_2$].
The second author would like to thank the hospitality of Osnabr\"uck
Universit\"at, SUNY \@ Stony Brook and  Stanford University during the
elaboration of portions of this work.

\subheading{\S 2. Spaces of complex cycles}  
For expository purposes we quickly review some known results for 
groups of algebraic cycles over $\bbc$.  The reader is referred to [\L$_2$] 
for an enlarged exposition.  Let $V$ be a  finite-dimensional complex vector 
space. 
For integers $d, q \geq 0$, let
$\calc_d^q(\bbp(V))$  denote   the Chow variety of effective 
algebraic cycles of
codimension $q$ and degree $d$ in the projective space $\bbp(V)$.  
The disjoint
union  $\calc^q(\bbp(V)) = \coprod_d \calc_d^q(\bbp(V))$ is an 
abelian
topological monoid whose na\"{\i}ve group completion is denoted by
$\calz^q(\bbp(V))$.  

As usual let $K(G,n)$ denote the 
Eilenberg-MacLane space with 
$\pi_nK(G,n) \cong G$ and $\pi_mK(G,n) \cong 0$ for $m\neq  n$,
and for a graded abelian group
$G_* = \oplus_{j\geq 0} \ G_j$, let $K(G_*)$ denote  the weak product 
$K(G_*)\ = \ \prod_{j\geq 0} K(G_j,j).$

\Theorem {2.1} ([\L$_1$]) {\sl There is a canonical homotopy 
equivalence}
$$
\cz^q(\bbp(V)) \ \cong \ K(\bbz,0) \times K(\bbz,2) \times K(\bbz,4) 
\times
\dots \times K(\bbz,2q) \tag{2.1}
$$
The canonical aspect of this splitting  is discussed in
Appendix A.

\Theorem {2.2} ([\LM])  {\sl The algebraic  join determines a 
continuous 
biadditive pairing
$$
\# : \cz^q(\bbp(V))\wedge \cz^{q'}(\bbp(V'))\ \arr\ 
\cz^{q+q'}(\bbp(V\oplus V'))
\tag{2.2}
$$
which, with respect to the splitting (2.1), represents the cup
product. }

\Theorem {2.3} ([\FM]) {\sl  Under the join pairing (2.2)  the 
homotopy groups
of the limiting space $\cz^{\infty}$ form a graded ring isomorphic to 
a
polynomial ring in one variable 
$$
\pi_* \cz^{\infty} \ \cong\ \bbz [s] 
\tag{2.3}
$$
where $s \in \pi_2\cz^{\infty}$ is the generator.}

\medskip

If one considers $\ \bbz [s]\ $  as a graded ring with one generator in
degree two, then the quotient $\ \bbz [s]/ (s^{q+1})\ $ has a natural
structure of graded abelian group. Using the terminology established
above, Theorems 2.1 and 2.3 can be reformulated by saying that one has
canonical homotopy equivalences 
$\cz^q(\bbp(V)) \ \cong \ K\left(\ \bbz [s]/(s^{q+1})\ \right) $ 
and 
$\cz^{\infty} \ \cong\ K\left( \bbz[s] \right)$. Furthermore,
the latter equivalence induces, under the join pairing, a
ring isomorphism $\pi_* \cz^{\infty} \ \cong\ \bbz [s]$.

Let $G^q(\bbp(V))= \cc^q_1(\bbp(V))$ denote the Grassmannian of 
codimension-$q$ planes in $\bbp(V)$, and let $\cz^q(\bbp(V))(1)$
denote the connected component of $\cz^q(\bbp(V))$ consisting of all
(not necessarily effective) algebraic cycles of degree $1$.

\Theorem {2.4} ([\LM]) {\sl  Under the splitting (2.1) the inclusion
$$
G^q(\bbp(V)) \hookrightarrow \cz^q(\bbp(V))(1)
\tag{2.4}
$$
represents the total Chern
class of the tautological $q$-plane bundle $\xi^q_{\bbc}$ over 
$G^q(\bbp(V))$.  
Passing to a limit as dim$(V) \to \infty$ gives a mapping
$$
BU_q \ \arr\ \cz^q(\bbp^{\infty})(1) \cong 
1 \ \times\ \prod_{i=1}^{q} K(\bbz, 2i)  
$$
which classifies the total Chern class of the universal 
$q$-plane bundle $\xi^q_{\bbc}$ over $BU_q$. Taking the limit as 
$q\to\infty$ gives 
a mapping
$$
BU \ \arr\ \cz^{\infty}(1) \cong 
1\ \times \ \prod_{i=1}^{\infty} K(\bbz, 2i) \ \equdef\ K(\bbz, 2*)
\tag{2.5}
$$
which classifies the total Chern class map from $K$-theory to even 
cohomology.}
\medskip

This natural presentation of the total Chern class map comes equipped
with the following remarkable property. 

\Theorem {2.5} ([\BLLMM]) {\sl   The join pairing on $K(\bbz, 2*)$
enhances to  an infinite loop space structure so that with respect to 
Bott's infinite loop structure on $BU$ the map (2.5) is an infinite
loop map.}

\subheading{\S 3.\ \ Spaces of real cycles}
A {\bf Real structure}  on a complex vector space $V$ is a
$\bbc$-antilinear map $\rho : V \rightarrow V$ such that
$\rho^2 = 1.$  A {\bf Real vector space}  is  a pair
$(V,\rho)$ consisting of a complex   vector space $V$
and a Real structure $\rho$.  Any such space is equivalent to
$(\bbc^n, \rho_0)$ where $\rho_0$ is complex conjugation.

A Real structure $\rho$ on $V$ induces
an anti-holomorphic $\bbz_2$-action on the
complex projective space $\bbp(V)$ which in turn induces an
anti-holomorphic $\bbz_2$-action on the Chow varieties
$\chv{q}{d}{\bbp(V)}$.
This produces an automorphism  
$$
\rho : \cyc{q}{\bbp(V)} \rightarrow \cyc{q}{\bbp(V)}.
\tag{3.1}
$$
of the topological group of all codimension-$q$ cycles on $\bbp(V)$.

\Def{3.1}
By the  {\bf Real algebraic cycles} of codimension $q$ on $\bbp(V)$ 
we mean
the subgroup  $\crl{q}{\bbp(V)}$ of cycles fixed by the involution 
$\rho$.
The closed subgroup of Galois sums  
$$
\cyav{q}{\bbp(V)}= \{ c + \rho c \; | \; c\in \cyc{q}{\bbp(V)}\}
$$
is called the group of {\bf averaged cycles}, and the quotient
$$
\crd{q}{\bbp(V)} = \crl{q}{\bbp(V)} /    \cyav{q}{\bbp(V)}
$$
is called the group of {\bf reduced Real algebraic cycles}.

\medskip

We have adopted  the standard definition of   
real algebraic cycles    as those which are fixed by the Galois group 
Gal$(\bbc/\bbr)$.  Note  that the group of
reduced   cycles is the topological vector space over $\bbz_2$
freely generated by the irreducible Real  subvarieties
of $\bbp(V)$.

Fix a Real vector space $(V,\rho)$ and let  $x_0 = [0:\cdots:0:1]
\in  \bbp(V \oplus \bbc)$.  Given an irreducible   algebraic 
subvariety  $Z
\subset \bbp(V)$ we define its {\bf algebraic suspension}   
 $\csus Z    = x_0\# Z \subset \bbp(V\oplus \bbc)$  to be the union 
 of all projective lines joining  $Z$ to $x_0$.
Algebraic suspension extends linearly to a $\bbz_2$ equivariant
 continuous homomorphism 
$$
\csus : \cyc{q}{\bbp(V)} \rightarrow \cyc{q}{\bbp(V\oplus \bbc)}.
\tag{3.2}
$$

The Algebraic Suspension Theorem  [\L$_1$] states that
(3.2) is a homotopy equivalence. When $V$ is a Real vector space,  T. 
K. Lam
showed that (3.2) is an {\sl equivariant}
homotopy equivalence.  
(See [\LLM$_2$] for considerable generalizations.)  In particular we 
have the
following.

\Theorem{3.2} \ ([\Lam]) \ \ {\sl  Algebraic suspension
induces  homotopy equivalences:
$$\aligned
\csus : \; \crl{q}{\bbp(V)}\  \xrightarrow{\cong}\ 
                &\crl{q}{\bbp(V\oplus \bbc)},  \qquad
\csus : \; \cyav{q}{\bbp(V)}\   \xrightarrow{\cong}\ 
                \cyav{q}{\bbp(V\oplus \bbc)},  \\ &\ \\
 \text{and}\qquad &\csus : \; \crd{q}{\bbp(V)}\  
\xrightarrow{\cong}\ 
                \crd{q}{\bbp(V \oplus \bbc)}.
\endaligned
$$
for all $q<\dim(V)$.}
\bigskip

This result shows that the homotopy types of the
topological groups $\crl{q}{\bbp(V)}$, $\cyav{q}{\bbp(V)}$ and 
$\crd{q}{\bbp(V)}$ depend  only on $q$, and so we can drop the 
reference to $V$. 
Our first theorems compute these homotopy types.  

\Theorem {3.3}\  {\sl
There is a canonical homotopy equivalence
$$
\cz_{\bbr}^q\ \cong \ \prod_{n=0}^q \prod_{k=0}^{n}
K(\I{n}{k}, n+k)
$$
where
$$\I{n}{k} = 
            \cases
            0 &,  \text{ if $k$ is odd or $k>n$};\\
            \bbz &, \text{ if $k=n$ and $k$ is even; } \\
            \bbz_2 &, \text{ if $k < n$ and  $k$ is even}.
             \endcases
$$
}

\Theorem{3.4}\ {\sl
There is a canonical homotopy equivalence
$$
\cz^q_{av}
\ \cong \ \prod_{n=0}^q \prod_{k=0}^{n}
K\left( \Iav{n}{k},\ n+k \right)
$$
where $ \Iav{0}{0} = 2\bbz$ and for $n+k>0$
$$\Iav{n}{k} \ =\ 
\tilde{H}_{k-1}(\pr^{n-1};\bbz)
            \ = \ \cases
            0 &, \text{ if $k$ is odd or $k>n$}; \\
            \bbz &, \text{ if $k=n$ and $k\geq 2$ is even; } \\
            \bbz_2 &, \text{ if $k< n$  and  $k\geq 2$ is even}.
             \endcases
$$
The homomorphism on homotopy groups induced by the inclusion
$\cz^q_{av} \subset \cz^q_{\bbr}$ is injective, and with respect 
to the splittings above, it maps $\Iav{n}{k}$ to $\I{n}{k}$
in the obvious way. (This explains the $2\bbz$ in $I_{0,0}$.)
}

 \medskip

\Theorem {3.5}\ ( [\Lam]) \ \ {\sl
 There is a canonical homotopy equivalence}
 $$
\widetilde{\cz}^q_{\bbr} \ \cong\ K(\bbz_2,0)\times 
K(\bbz_2,1)\times
 K(\bbz_2,2)\times\dots \times  K(\bbz_2,q).
\tag{3.1}
$$
The proofs of these results are given in \S 8.
Useful diagrams of the graded groups ${\ci}^{av}_{*,*}$ and 
$\ci_{*,*}$
are given in \S 9.

\subheading{\S 4.\ The ring structure}
The homotopy groups
$$
\pi_* {\cz}^q_{\bbr}\  = \
\bigoplus_{0\leq k\leq n\leq q} \I {n}{k}
\tag{4.1}
$$
are vastly simplified conceptually if one takes into account 
their {\bf multiplicative structure}.  The algebraic join pairing (2.2)
restricts to a pairing
$$
\# : {\cz}^q_{\bbr}\wedge {\cz}^{q'}_{\bbr} \ \arr\ {\cz}^{q+q'}_{\bbr}
$$ 
which gives $\pi_* {\cz}^{\infty}_{\bbr}$ the structure of a 
commutative ring. Since the join of an averaged cycle with a fixed
cycle is again an  averaged cycle, the subgroup $\pi_*
{\cz}^{\infty}_{av}$ is an {\sl ideal} in this  ring.

In \S 9 we will prove the following result.

\Theorem {4.1}  {\sl 
There is a ring isomorphism
$$
\pi_* {\cz}^{\infty}_{\bbr} \ \cong\ \bbz[x,y]/(2y)
\tag{4.2}
$$
where $x$ corresponds to the generator of $\pi_4 
{\cz}^{\infty}_{\bbr}  \cong
\bbz$ and $y$ corresponds to the generator of $\pi_1 
{\cz}^{\infty}_{\bbr}  \cong
\bbz_2$, and where $(2y)$ denotes the principal ideal generated by 
$2y$  in the
polynomial ring $\bbz[x,y]$.  Under this isomorphism  the ideal 
$\pi_*{\cz}^{\infty}_{av} \subset \pi_*{\cz}^{\infty}_{\bbr}$ 
corresponds to the
ideal 
$$
\pi_*{\cz}^{\infty}_{av}\ \cong\ (2,x)
$$
generated by 2 and $x$.
Furthermore, with respect to the isomorphisms (4.1) and (4.2), we 
have
$$
{I}_{2m+\ell,2m} \text{\ \  is the cyclic subgroup 
generated by \ \ }
x^m y^{\ell}.  
$$
}

\Cor {4.2}\ {\sl The algebraic join induces a ring structure on
$\pi_*\widetilde{\cz}^{\infty}_{\bbr}$. There is a canonical 
ring isomorphism
$$
\pi_*\widetilde{\cz}^{\infty}_{\bbr}\ \cong \ \bbz_2 [y]
$$
where $y$ in the generator of $\pi_1\widetilde{\cz}^{\infty}_{\bbr}
=\bbz_2$.}

\Remark {4.3}\ Consider the polynomial ring $\bbz[x,y]$ on the variables
$x$ and $y$, of degrees $4$ and $1$, respectively. Given a non-negative
integer $q$, define the ideal  
$$\calj_q = ( \ 2y, \ \{ x^my^j\ : \ 2m+j = q+1\ \} ) \subset \bbz[x,y],$$ 
and denote $\calj_\infty = ( 2y )$. Each  quotient ring 
$\calr^q_*=\bbz[x,y]/\calj_q$, $q=0,\ldots,\infty$, has the  natural
structure of a graded abelian group. 

Using this notation, Theorem 3.3 can be rephrased by saying that there
is a canonical equivalence 
$$\calz^q_\bbr \cong K(\calr^q_*).$$
Under this equivalence the direct summand $I_{2q+j,2q}$ of the
$(4q+j)$-th homotopy group of  $\calz^q_\bbr$ is precisely the
subgroup of $\calr^q_*$ generated by $x^qy^j$. One can rephrase
Theorem 3.4 and Corollary 3.5 in a similar fashion.

We also prove that there are canonical equivalences 
$$
{\cz}^{\infty}_{\bbr} \ \cong\ K\left( \ \calr^\infty_*  \ \right), 
\quad
\widetilde{\cz}^{\infty}_{\bbr}\ \cong \ K\left( \ \bbz_2 [y]\ \right) 
\quad \text{ and } \quad 
{\cz}^{\infty}_{av}\cong K(\ I^{av}\ ).
$$
Here $I^{av}$ is the ideal $I^{av} = (2,x) \subset \calr^\infty_*$.
Furthermore,  these homotopy equivalences induce the ring isomorphisms
presented in Theorems~4.1 and Corollary 4.2.
\medskip

\subheading{\S 5.\ Extending functors from $K$-theory}
We shall now show
that certain basic functors in classical representation
theory carry over to algebraic cycles.  
This remarkable fact together with   [\LM] and the results
of \S 3  leads to a new proof of the basic relationships among
characteristic classes.

Before beginning we set some notation.  For all $k\geq 0$ let
$$
\iota_{2k} \in H^{2k}(K(\bbz,\, 2k); \, \bbz)\cong \bbz
\qquad\text{and}\qquad
\widetilde{\iota}_k \in H^{k}(K(\bbz_2,\, k); \, \bbz_2)\cong \bbz_2
$$
denote the {\sl fundamental classes} (i.e., the canonical generators).
Let $c_k$, $w_k$, and $p_k$ denote respectively the k$^{\text{th}}$ 
Chern,
Stiefel-Whitney, and Pontrjagin classes.

\medskip

\subheading{Complexification}  Consider a Real vector space 
$(V,\rho)$
and the map $(V,\rho) \to V$ which forgets the Real structure.  
Associated to 
this is the homomorphism $\crl q {\bbp(V)} \hookrightarrow 
\cz^q_{\bbc}(\bbp(V))$
which simply includes the subgroup  fixed by   $\rho$ into the group 
of all 
cycles. Restricting to linear  cycles  gives
a commutative diagram
$$\CD
G^q_{\bbr}(\bbp(V))  @>>>     G^q_{\bbc}(\bbp(V)) \\
@V{P}VV     @VV{c}V     \\
\cz^q_{\bbr}(\bbp(V))  @>>>     \cz^q_{\bbc}(\bbp(V))
\endCD
\tag{5.1}
$$
where 
$$
G^q_{\bbr}(\bbp(V)) \ =\ \{\ell \in G^q_{\bbc}(\bbp(V))\, :\, 
\rho(\ell) =
\ell\}
$$
is the Grassmannian of real subspaces of codimension-$q$ in
$
V_{\bbr}  = \{v\in V\,:\, \rho(v) = v\}.
$

Recall from 2.4 that under the canonical identification 
$ 
\cz^q_{\bbc}(\bbp(V)) \ \cong \ \prod_{k=0}^q K(\bbz,\, 2k)
$ 
the map $c$ in (5.1) classifies the total Chern class of the 
tautological $q$-plane bundle                                \linebreak
$\xi^q_{\bbc} \arr G^q_{\bbc}(\bbp(V))$,
i.e.,  $$
c^*(\iota_{2k}) = c_k(\xi^q_{\bbc}) \qquad\text{for}\ k=0, ..., q.
$$

Consider   the composition  
$$
w = \pi\circ P:G^q_{\bbr}(\bbp(V))\arr \crdh {q}(\bbp(V))
\tag{5.2}
$$ 
where $\pi : \cz^q_{\bbr}(\bbp(V)) \to \crdh {q}(\bbp(V)) = 
\cz^q_{\bbr}(\bbp(V))/\cz^q(\bbp(V))^{\text{av}}$ is the projection.
It is a result of Lam [\Lam] that under the canonical identification 
$ 
\crdh {q}(\bbp(V)) \ \cong \ \prod_{k=0}^q K(\bbz_2,\, k),
$ 
the map $w$   classifies the total Stiefel-Whitney class of the 
tautological real $q$-plane bundle $\xi^q_{\bbr} \arr 
G^q_{\bbr}(\bbp(V))$, i.e.,
$$
w^*(\widetilde{\iota}_k) = w_k(\xi^q_{\bbr}) \qquad\text{for}\ k=0, 
..., q.
$$

\medskip

 We now set $V = \bbc^n$ and take the colimit of the spaces in (5.1) 
and (5.2)
as $n\to \infty$.  This gives a diagram
$$
\CD
BO_q   @>{\cx}>>    BU_q  \\
@V{P}VV     @VV{c}V   \\
\cz^q_{\bbr} @>{\cxx}>>  \cz^q_{\bbc}\\
@V{\pi}VV     @.   \\
\crdh {q} @. \
\endCD
\tag{5.3}
$$
where $\cz^q_{\bbc} \equiv 
\underset{n\to\infty}\to{\lim}\cz^q_{\bbc}(\bbp(\bbc^n))$,  etc..
By using (3.3), this can be canonically rewritten  as
$$
\CD
BO_q   @>{\cx}>>    BU_q  \\
@V{P}VV     @VV{c}V   \\
{\dsize \prod_{n=0}^{[q/2]}\prod_{i=1}^{q-2n}K(\bbz_2,4n+i)}  \times 
{\dsize\prod_{k=0}^{[q/2]}} K(\bbz, 4k)
 @>{\cxx}>> {\dsize \prod_{k=0}^{q}}K(\bbz,2k)\\
@V{\pi}VV     @.  \\
{\dsize\prod_{k=0}^{q}}K(\bbz_2,k) @. \
\endCD
\tag{5.4}
$$
The map $\cx$ on classifying spaces is the one induced by the 
inclusion
$O_q \subset U_q$ associated to the complexification of vector spaces
$V_{\bbr} \mapsto V_{\bbr}\otimes \bbc$.  

Consider the classes 
$$
\pont_{2k} \ \equdef\ \cxx^*\iota_{2k} \,\in\, 
H^{2k}(\cz^q_{\bbr};\,\bbz)
\qquad\text{and}\qquad
\whit_{k} \ \equdef\ \pi^*\widetilde{\iota}_{k} \,\in\,
H^{k}(\cz^q_{\bbr};\,\bbz_2).
$$
From these theorems and the commutativity of the diagrams above 
we see that
$$
(-1)^kP^*\pont_{4k} \ =\ p_k(\x^q_{\bbr})
\qquad\text{and}\qquad
P^*\whit_{k} \ =\ w_k(\x^q_{\bbr}).
\tag{5.5}
$$
In particular, $\pont_{4k}$ is not divisible and not torsion, whereas $j_{4k+2}$ has order $2$.  From the   factoring  
(5.10), (5.11) below we see that $\pont_{4k} = \iota_{4k} + \tau$ 
where $2\tau = 0$.

\vskip .3in

\subheading{The forgetful functor}  For a complex vector space $V$ 
one
constructs the conjugate space $\overline{V}$ by taking the same 
additive group
and defining a new scalar multiplication $\bullet$ by $t\bullet v 
\equiv
\overline{t} v$.  With this we can associate to $V$ a Real space
$([V]_{\bbr}, \rho)$ where
$$
[V]_{\bbr}\ =\ V \oplus \overline{V}
\qquad\text{and}\qquad
\rho(v,w) = (w,v).
$$
For any $q<\dim(V)$ we have a map
$$
\Phi : \cz^q_{\bbc}(\bbp(V))\ \arr\ \cz^{2q}_{\bbr}(\bbp(V\oplus 
\overline{V}))
$$
defined by
$$
\Phi(c) \ =\  c \# c
$$
where $\#$ is the complex join.
This construction gives  commutative diagrams
$$\CD
G^q_{\bbc}(\bbp(V))  @>{\Real}>>     G^{2q}_{\bbr}(\bbp([V]_{\bbr})) 
\\
@V{c}VV     @VV{P}V     \\
\cz^q_{\bbc}(\bbp(V))  @>{\Reall}>>     
\cz^{2q}_{\bbr}(\bbp([V]_{\bbr}))
\endCD
$$
which stabilize as above to commutative diagrams
$$
\CD
BU_q   @>{\Real}>>    BO_{2q}  \\
@V{c}VV     @VV{P}V   \\
\cz^q_{\bbc} @>{\Reall}>>  \cz^{2q}_{\bbr}.
\endCD
$$
Note that $\Real$ is the map induced by the standard inclusion
$U_q \subset O_{2q}$.

\vskip .3in
\subheading{Relations}  Consider the diagram
$$\CD
BU_q @>>>  BO_{2q} @>>>  BU_{2q} \\
@V{c}VV   @V{P}VV   @VV{c}V   \\
\prod_{j=0}^q K(\bbz,2j)  @>{\Phi}>> \cz_{\bbr}^{2q} @>{\Gamma}>>
\prod_{j=0}^{2q} K(\bbz,2j) 
\endCD
\tag{5.6}
$$
Note  that if $V$ has a real structure $\rho$, then under the 
isomorphism $I\oplus \rho:V \oplus \overline V \arr V \oplus V$, 
the map 
${\Phi} :\cz_{\bbc}^{q} @>>> \cz_{\bbr}^{2q}$ becomes 
${\Phi}(c) = c \# \rho_*(c)$.  It follows that 
$$
\Gamma \circ{\Phi}(c) = c \# \rho_*(c)
$$ 
for $c \in \cz_{\bbc}^{q}$.
We conclude the following.

\Prop{5.1} {\sl The  composition $\Gamma \circ{\Phi}$ satisfies
$$\aligned
(\Gamma \circ{\Phi})^* \iota_{2k} \ &= \ \sum_{i+j=k} (-1)^j\iota_{2i} 
\cup
\iota_{2j}  \\
&\ = \ 
\cases
2\sum_{j=0}^{m-1} (-1)^j \iota_{2j}\cup\iota_{2(2m-j)} + (-1)^m 
\iota_{2m}^2
&\text{if} \ \ k=2m \\
0 &\text{if} \ \ k=2m+1
\endcases
\endaligned
\tag{5.7}
$$
for all $k$.  }

\pf
By Theorem 2.2 the  join mapping 
$\# : \cz_{\bbc}^q\times \cz_{\bbc}^q \arr\cz_{\bbc}^{2q}$
has the characterizing property that
$
\#^* \iota_{2k}\ =\ \sum_{i+j=k}  \iota_{2i} \otimes \iota_{2j}.
$  
It is straightforward to verify that the map
$\rho: \cz_{\bbc}^q \to \cz_{\bbc}^q$, induced by the real structure 
$\rho$,
has the characterizing property that
$
\rho^* \iota_{2k} \ =\ (-1)^k \iota_{2k}.
$
Taking the composition
$\cz_{\bbc}^q @>{\Delta}>> \cz_{\bbc}^q\times \cz_{\bbc}^q 
@>{1\times \rho}>> \cz_{\bbc}^q\times \cz_{\bbc}^q 
@>{\#}>>\cz_{\bbc}^q$
and pulling back $\iota_{2k}$ gives the result.  \qed

\vskip .3in

 Similarly we have the  diagram
$$\CD
BO_q @>{\gamma}>>  BU_{q} @>{\phi}>>  BO_{2q} \\
@VVV   @VVV   @VVV   \\
\cz_{\bbr}^{q}  @>{\Gamma}>> \prod_{j=0}^q K(\bbz,2j)  @>{\Phi}>>
\cz_{\bbr}^{2q}
\endCD
\tag{5.8}
$$
and the relation
$$
\Phi\circ \Gamma (c) = c\# \rho_* c = c \# c,
$$
i.e., $\Phi\circ \Gamma $ is just the {\bf squaring map}. 
Thus $\Phi\circ \Gamma$ induces a map
$$
\widetilde{\Phi\circ \Gamma} : \widetilde{\cz}_{\bbr}^q
\arr\widetilde{\cz}_{\bbr}^{2q}  
$$
which is also the squaring map.  

\Prop{5.2} {\sl  Under the canonical equivalence 3.5, the
composition $\Phi\circ \Gamma$ satisfies 
$$
(\widetilde{\Phi\circ \Gamma})^* \widetilde{\iota}_k \ = \ 
\sum_{i+j=k} 
\widetilde{\iota}_i \cup\widetilde{\iota}_j  \ = \ 
\cases
\widetilde{\iota}_m^2  &\text{if} \ \ k=2m \\
0 &\text{if} \ \ k=2m+1
\endcases
\tag{5.9}
$$
for all $k$.}

\pf
Using the fact that the join $\# :\widetilde{\cz}_{\bbr}^q
\times \widetilde{\cz}_{\bbr}^q \arr \widetilde{\cz}_{\bbr}^{2q}$
classifies the cup product [\Lam], one proceeds as in the proof of 5.1. 
\qed
\vskip .3in

Note that the composition $BO_q \to \widetilde{\cz}_{\bbr}^q \to  
\widetilde{\cz}_{\bbr}^{2q}$ classifies the square of the total 
Stiefel-Whitney class $w(\xi_{\bbr}^q)^2 = \sum_{k=0}^q 
w_k(\xi_{\bbr}^q)^2$.

\vskip .3in

Notice that as $q$ increases the diagrams (5.8) are included in one 
another.
From [\L$_1$] and [\Lam] we know that 
$$
\cz^q_{\bbc}/\cz^{q-1}_{\bbc}\cong K(\bbz,2q)
\qquad\text{and}\qquad
\widetilde{\cz}^q_{\bbr}/\widetilde{\cz}^{q-1}_{\bbr}\cong 
K(\bbz_2,q).
$$
From this we obtain a diagram
$$
\CD
BO_q/BO_{q-1} @>>>  BU_q/BU_{q-1}  \\
@VVV   @VVV  \\
{\cz}^q_{\bbr}/{\cz}^{q-1}_{\bbr} @>{\Gamma_0}>>  K(\bbz,2q) \\
@V{\pi_0}VV @.\\
K(\bbz_2,q). @. \ 
\endCD
\tag{5.10}
$$
and from Theorem 3.3 we know that
$$\aligned
{\cz}^{2q_0}_{\bbr}/{\cz}^{2q_0-1}_{\bbr}\ &=\ 
K(\bbz, 4q_0) \times \prod_{i=0}^{q_0-1} K(\bbz_2, 2q_0+2i) 
\qquad\text{and}\\
{\cz}^{2q_0+1}_{\bbr}/{\cz}^{2q_0}_{\bbr} \ &=\ 
\prod_{i=0}^{q_0} K(\bbz_2, 2q_0+2i+1).
\endaligned\tag{5.11}$$
By Theorem 3.4 the map $\pi_0$ kills all factors with $i>0$.
However, $\Gamma_0$ could represent non-trivial cohomology 
operations on
$K(\bbz_2, 2*)$.

\medskip
\subheading{Question 5.3}  What is the cohomology class 
$\Gamma_0^*( \iota_{2q})$?
\medskip

We can say something about its structure.

\Prop {5.4}  {\sl  The class $\G_0^*(\iota_{2q})$ is ``primitive'',  i.e.,   
it is of the form   
$$
\G_0^*(\iota_{2q})\ =\ 
\cases
\iota_{2q}\otimes 1 \otimes\dots\otimes1 +\sum_{i=0}^{q_0-1}
1 \otimes\dots\otimes1\otimes x_i\otimes1 \otimes\dots\otimes1 
& \text{when $q = 2q_0$}\\
\sum_{i=0}^{q_0}1 \otimes\dots\otimes1\otimes x_i\otimes1 
\otimes\dots\otimes1 
& \text{when $q = 2q_0+1$}
\endcases
$$
where $x_i \in H^{2q}(K(\bbz_2, q+2i);\,\bbz)$.
}
\pf
Suppose $q = 2q_0+1$ and set $K_i = K(\bbz_2, q+2i)$. 
Define $\calk(r) \equdef \prod_{j=q_0-r}^{q_0} K_j$.
Since $\pi_k(\calk(r)) = 0$ if $k< q+2(q_0-r)$ and
$\pi_{q+2(q_0-r)}(\calk(r)) = \bbz_2$, it follows from Hurewicz
theorem and the universal  coefficients theorem that 
$H^{\alpha}(\calk(r)) = 0$ if 
$$\alpha \leq q+2(q_0-r).
\tag{5.12}$$

We now prove, using induction, that if $\alpha \leq 2q$ then
$$
H^{\alpha}(\calk(r)) \cong 
\left\{
H^*(K_{q_0-r})\otimes \dots\otimes H^*(K_{q_0})\right\}^{\alpha},
\tag{5.13}
$$
for all $0\leq r\leq q_0$.
 
For $r=0$ the assertion is evident. Assume (5.13) holds for
all $0\leq r <m \leq q_0$, and note that
$\calk(m) = K_{q_0-m} \times \calk(m-1)$. Now, recall the 
Kunneth exact sequence
$$
\aligned
0&
\arr\left\{ H^{*}(K_{q_0-m})\otimes  H^{*}(\calk(m-1))\right\}^{\alpha} \arr
H^{\alpha}(K_{q_0-m}\times \calk(m-1)) \\
&\arr \left\{ H^{*}(K_{q_0-m})* H^{*}(\calk(m-1))\right\}^{\alpha +1} \arr 0 
\endaligned
\tag{5.14}
$$
where $*$ denotes the tor-product. Since 
$
H^{m_i}(K_i;\bbz)
$ is free or zero unless
$$
m_i \geq q+2i,
\tag{5.15}
$$
it follows from (5.15) and (5.12) that 
$ H^{a}(K_{q_0-m})* H^{b}(\calk(m-1)) = 0$
unless $a+b \geq q+2(q_0-m) + q +2\{ q_0- (m-1)\} +1 =
2q +4(q_0-m)+3 \geq 2q+3$. It follows that the tor-term in (5.14) is
zero if $\alpha \leq 2q$. 

Since 
$$\left\{ H^{*}(K_{q_0-m})\otimes  H^{*}(\calk(m-1))\right\}^{2q}=
\oplus_{a+b=2q}H^a(K_{q_0-m})\otimes  H^{b}(\calk(m-1)),
$$
and $b \leq 2q$, the induction hypothesis completes the proof of (5.13).

Observe now that the tensor product 
$H^{m_0}(K_0)\otimes\dots\otimes  H^{m_{q_0}}(K_{q_0})$
is non-zero only when each $m_i$ is either 0 or satisfies (5.15).
If $m_i$ and $m_j$ satisfy (5.15) with $i\neq j$, then as above 
$m_i+m_j 
\geq 2q + 2 > 2q$.  This proves the Proposition when $q$ is odd.

The case when $q$ is even is similar. One must only see that the 
contribution
from $H^{2q}(K(\bbz, 2q);\,\bbz) = \bbz\cdot \iota_{2q}$ is exactly
$\iota_{2q}$.  This follows from the fact that $P^*\G^*\iota_{2q}$
is, up to sign,  the $q^{\text{th}}$ Pontrjagin class of the universal 
bundle over BO$_{2q}$, which is indivisible. This is seen from  the
commutativity of (5.6) and the fact that $c$ represents the total 
Chern class.
 \qed

\vskip .3in

Consider now the composition  $P$ given by
$$
  BO_q \ar{Q} \ BO_q/BO_{q-1} \ \ar{P_0} \ 
\cz_{\bbr}^q/\cz_{\bbr}^{q-1}
$$
where $Q$ is the quotient and $P_0$ comes from (5.10). The image of 
$Q^*$
in $H^*(BO_q;\,\bbz_2)  \cong \bbz_2[w_1,\dots,w_q]$ is the ideal 
$(w_q)$
generated by the $q^{\text{th}}$ Stiefel-Whitney class $w_q$.  
Consider the
canonical product structure (5.11) and let 
$K(\bbz_2;x^iy^{q-2i})$ denote the Eilenberg-MacLane space
$K(\bbz_2;q+{2i})$ corresponding to the monomial
$x^iy^{q-2i}$ under the identification
$\calz_\bbr \cong K\left( \bbz[x,y]/(2y) \right).$
With this notation, let
$\widetilde{\iota}_{q,2i}$
denote the fundamental class in 
$H^{q+2i}(K(\bbz_2, x^iy^{q-2i})\, ;\, \bbz_2)$ 
pulled back to the product.  Note that $\tilde{\iota}_{q,2i}$ 
is the Kronecker dual to the class ${\theta}_{q,2i} $ introduced in (9.11). 
Then we see that
$$
P^*\widetilde{\iota}_{q,2i} \ =\  F_{q,2i}(w_1,\dots,w_{q-1})\cdot 
w_q
\tag{5.16}
$$
where $F_{q,2i}(\xi_1,\dots,\xi_{q-1}) \in \bbz_2[\xi_1,\dots,\xi_{q-
1}]$ is a
homogeneous polynomial of weighted degree $2i$, i.e., 
$F_{q,2i}(t\xi_1,t^2\xi_2,\dots,t^{q-1}\xi_{q-1}) = t^{2i} F_{q,2i}(\xi)$.
These polynomials determine $P$ up to homotopy.

\subheading{\S 6.  Equivariant infinite loop space structures and
$KR$-theory}  
In this section we  shall show that our spaces of complex algebraic 
cycles have the structure of an 
{\bf equivariant $E_{\infty}$-ring spaces}  
(cf. [\LMS]), under the $\bbz_2$ action induced by complex
conjugation. The principle is the same as in [\LLM$_1$], where the
ruled join of cycles induces the infinite loop structure. However,
here we obtain  $RO(\bbz_2)$-graded cohomology theories, as opposed to
$R(\bbz_2)$-graded ones.

Furthermore, we show that one obtains two canonical equivariant
infinite loop spaces from these constructions. The first one comes
from delooping the additive structure, which yields an equivariant
ring spectrum. The second one come from delooping the multiplicative
units of the original ring space. This yields an equivariant spectrum
which is directly related to characteristic classes in Atiyah's
$KR$-theory.

It follows from these constructions that the space of real cycles 
$\cz_{\bbr}^{\infty}$ is also an $E_{\infty}$-ring space and that most of
the maps introduced in previous sections are maps of $E_\infty$-ring
spaces. Our arguments involve P. May's use of equivariant
$\ci_*$-functors and make  extensive use of the constructions in
[\LLM$_1$].  We shall briefly introduce the  concepts  but refer the
reader to [\LLM$_1$] for many details. 

Consider $\bbc^\infty$ as a direct sum 
$\bbr^{\infty} \!  \oplus i\bbr^\infty$ with its
standard orthogonal inner product, and where $\bbz_2$ acts by complex
conjugation.  Then $\bbc^\infty$ contains
infinitely many copies of each irreducible real representation of
$\bbz_2$, in other words,  in the terminology of [\LMS] it is a {\it
complete $\bbz_2$-universe}. It will be fixed throughout this
discussion. 

In general, suppose $G$ is a finite group and let $\cu$ be a fixed
$G$-universe. Recall that an   {\bf equivariant} infinite loop space
$X$, indexed on  $\cu$, is a based $G$-space for which there is
collection of $G$-spaces
$\{ X(V) \ | \ \text{ $V \subset \cu $ is a $G$-submodule} \}$  
together with  $G$-equivariant homeomorphisms
$X\cong \Omega^V X(V)$. Here $ \Omega^V X(V)$
denotes the space of based maps $F(S^V,X(V))$ from the one-point
compactification $S^V$ of $V$ to $X(V)$, and $ \Omega^V X(V) $ is
equipped with its natural structure of $G$-space. 
The structural homeomorphisms are coherent in the sense that, if for a given
submodule $W   \subset V$ one denotes by $V-W$ the orthogonal
complement of $W$ in $V$,  then there are compatible
$G$-homeomorphisms $X(W) \cong \Omega^{V-W} X(V).$ In general, to give
a $G$-space $X'$  an equivariant infinite loop space structure is to
provide a $G$-homotopy equivalence between $X'$ and an equivariant
infinite loop space $X$.

\Remark{6.1}  If $\{ X(V) \ | \ V\subset \cu \}$ is the collection of
equivariant ``deloopings'' of the equivariant infinite loop space $X$,
let $X(n)$ denote $X(\bbr^n)$ for the trivial $G$-module
$\bbr^n$. Then for any subgroup $H\leq G$, the fixed point set
$X^H$ has the structure of a (non-equivariant) infinite loop space,
since the $G$-homeomorphism $X \cong \Omega^n X(n)$ gives a
homeomorphism of fixed point sets $X^H \cong \left(  \Omega^n X(n)
\right)^H =  \Omega^n\left( X(n)^H \right).$ Furthermore, if $H\leq K
\leq G$ are subgroups then, under the structure defined above, the
inclusion $X^K \subset X^H$ is obviously a map of (non-equivariant)
infinite loop spaces.
\medskip

In order to show that a $G$-space $X$ has the structure of an
equivariant infinite loop space, we use the machinery 
developed in  [\CW]. In this formulation, one considers the
category of {\bf $G\cl(\cu)$-spaces}, whose objects are $G$-spaces on
which there is an action of the equivariant linear isometries operad
$G\cl(\cu)$  (cf. [\May$_3$, pp 10 ff], [\CW]), 
and where a {\bf map of $G\cl(\cu)$-spaces} is a $G$-map which
commutes with the action of $G\cl(\cu)$.  The next result is a
formulation of the main results from [\CW], suitable for our purposes.

\Theorem{6.2} ([\CW]) {\sl 
Let $\cu$ be a complete $G$-universe and let $X$ be a
$G\cl(\cu)$-space which is $G$-group-complete. In other words,
for each subgroup $K\leq G$ the induced $\Cal H$-space structure
makes $\pi_0(X^K)$ a group. Then $X$ has
an equivariant infinite loop space structure. This structure is
natural in the sense that any map of  $G$-group-complete
$G\cl(\cu)$-spaces induces an equivariant infinite loop map.    
}
\medskip

From now on, we restrict ourselves to the case where $G= \bbz_2$ and
fix the $\bbz_2$ universe $\cu = \bbc^{\infty}$ described above. For
simplicity we shatll write $\cl$ instead of $\bbz_2\cl(\cu)$, and we
shall avoid mentioning the universe in most instances.

Following Atiyah's terminology [\A], define a 
{\bf Real topological space} to be a pair $(X,\rho)$ where 
$X$ is a space and $\rho:X\to X$ is an involution. In other 
words, $X$ is a $\bbz_2$-space. A {\bf Real mapping} $f:X\to Y$
between Real spaces is one which  commutes with the involutions
(a $\bbz_2$-equivariant map).  We denote by  $\z2t$  the
category  of compactly generated, based Hausdorff  Real
topological spaces, with base-point fixed by the action. The
morphism spaces in $\z2t$  are given the usual topology in the compactly
generated category, and have the natural $\bbz_2$-action on them.  
\medskip

A natural way of constructing actions of the equivariant linear
isometries operad $\cl$ uses the following notions. Let $\bbz_2\ci_*$
the subcategory of the category of finite dimensional hermitian
$\bbz_2$-modules and $\bbz_2$-module morphisms, whose morphisms are
also linear isometries.

\Def{6.3}  A {\bf  $\bbz_2 \ci_*$-space} (or $\bbz_2 \ci_*$-{\bf functor}) 
$(T,\omega)$ is a continuous
covariant functor $T:\bbz_2\ci_* \arr \z2t$ together with 
a (coherently) commutative, associative and continuous natural
transformation  
$\omega:  T\times T  \arr T \circ \oplus$ such that
\roster
\item   If $x\in TV$ and if $1 \in T\{0\}$ is the basepoint, then
$$
\omega(x,1) = x \in T(V\oplus\{0\}) = TV,
$$
\item  If $V = V'\oplus V''$, then the map $TV'\arr TV$ given by 
$x\mapsto
\omega(x,1)$ is a homeomorphism onto a closed subset;
\item Each sum map $\omega : T(V) \times T(W) \to T(V\oplus W)$ is a
$G$-map;
\item Each evaluation map $e: \bbz_2\ci_*(V,W) \times T(V) \to T(W)$
is a $G$-map.
\endroster
\medskip

The following result is a direct consequence of the techniques in 
[\May$_3$]. See the discussion in [LLM$_1$, \S 2].

\Theorem{6.4} \  \ \ {\sl If $(T,\omega)$ is an
$\bbz_2\ci_*$-space, then
$$
T(\bbc^{\infty}) \ = \ \lim_{V\subset \bbc^{\infty}} \, T(V),
$$
where the limit is taken over finite-dimensional $\bbz_2$-submodules
of  $\ \bbc^{\infty}$, is an $\cl$-space.  Any map
$\Phi: (T,\omega) \arr (T',\omega')$ 
of  $\ci_*$-spaces, induces a mapping  $\Phi:T(\bbc^{\infty})
\arr T'(\bbc^{\infty})$ of $\cl$-spaces.}
\medskip

A given $V \in \bbz_2\ci_*$ can be written as $V=\bbr^n \oplus
\sigma\otimes \bbr^m$, where $\bbr^k$ denotes a trivial representation
of rank $k$ and $\sigma$ is the sign representation of $\bbz_2$. In
particular, if one denotes by $V_{\bbr}$ the underlying real vector
space of $V$, then the sum $V \oplus \sigma\otimes V$ is canonically
isomorphic to $V_{\bbc} \equdef V_{\bbr}\otimes \bbc$ as a
$\bbz_2$-module, where the action on the latter is given by complex
conjugation. Given such $V$, we denote its real dimension by $v=n+m$,
and for any map $f: V \to W$ we denote by $f_{\bbc}$ its natural
extension  to the complexified vector spaces.

\Ex{6.5} {\bf (The Grassmann functor)} \ \ Given $V\in \ \bbz_2\ci_*$ of 
dimension $v$, let $T_G(V) = G^v(V_{\bbc}\oplus
V_{\bbc})=G^v(\bbp(V_{\bbc}\oplus V_{\bbc}))$ be the  Grassmannian of
codimension-$v$ complex planes in $V_{\bbc}\oplus V_{\bbc}$,
with distinguished point $1_G = V_{\bbc}\oplus \{0\}$.  To a linear isometric 
embedding $f:V \to W$ we define $T_Gf: T_GV\to T_GW$ on a plane 
$P\subset V_{\bbc}\oplus V_{\bbc}$ 
by 
$T_Gf(P) = ((f_{\bbc}V_{\bbc})^{\perp}\oplus \{0\}) \oplus (f_{\bbc}\oplus f_{\bbc})(P)$. 
The natural transformation $\omega_G : T_G\times T_G \arr
T_G \circ \oplus$ is given by the direct sum, i.e., for $P \in T_GV$ 
and
$P' \in T_GV'$, we define 
$\omega_G(P,P') = \tau_{*}(P\oplus P')$ 
where  
$\tau : V\oplus V\oplus V'\oplus V' \arr V\oplus V'\oplus V\oplus V'$
is the isometry interchanging the  middle factors.
This is an $\bbz_2\ci_*$-functor, and  
$$
T_G(\bbc^{\infty}) \ = \ \text{\bf  BU}
$$
is then a $\bbz_2$-equivariant $\cl$-space,  and hence it is an
equivariant infinite loop space; cf. Theorem 6.2. 
According to Remark 6.1, if $\{ 0 \}$ denotes the trivial subgroup of
$\bbz_2$, then both fixed point sets 
$\text{\bf BU} = \text{\bf BU}^{\{ 0 \} } $ and  
$\text{\bf B0} = \text{\bf BU}^{\bbz_2}$ inherit infinite loop space
structures which makes the  canonical ``complexification'' inclusion
$\text{\bf BO}\hookrightarrow \text{\bf BU} $ a map of infinite loop spaces.  These are the
standard Bott infinite loop space structures.  (See [\May$_3$,
pg.16].)  
\medskip

This equivariant structure on $\text{\bf BU}$ classifies an
$RO(\bbz_2)$-graded equivariant cohomology theory which we now
recall. 

\Def{6.6} \ Let $(X,\rho)$ be a Real space.  
A {\bf Real vector bundle} over $(X,\rho)$ is a Real space 
$(E,{\rho_E})$ where $\pi : E\to X$ is a  complex vector bundle,
${\rho_E}$ is a complex {\sl anti-linear} bundle map, and $\pi$ is a Real map,
i.e., $\pi{\rho_E} =\rho \pi$.  
\medskip

A Real projective variety  with its complex conjugation  involution
gives a Real  space. Important examples are the  Grassmannians
$G^q(\bbc^n)$ and the  Chow varieties.  The universal $q$-plane bundle
$\xi^q$ over  $G^q(\bbc^n)$ is a Real bundle.

\Prop{6.7} {\sl  Let $(X,\rho)$ be a  Real space which is compact 
and Hausdorff. Then the association $f \mapsto f^*\xi^q$ gives an 
equivalence of functors: 
$$
[X, G^q(\bbc^{\infty})]_{\bbr}\  @>{\cong}>> \ 
{\roman Vect}^q_{\bbr}(X), 
$$
from homotopy classes of Real mappings
$X\to G^q(\bbc^{\infty})$ to the set ${\roman Vect}^q_{\bbr}(X)$  of
equivalence classes of Real $q$-dimensional vector bundles over 
$X$.}

\pf
One can carry through the standard proof (cf. [\MS]).  The only point
to establish is that a Real bundle is locally trivial in the category
of Real bundles. This is shown for example in [\A]. \qed
\medskip 

It follows that the limiting Real space $G^{\infty} \cong
G^{\infty}( \bbc^{\infty}\oplus \bbc^{\infty} ) = \text{\bf BU}$
classifies  Atiyah's $KR$-theory ([\A]), 
and hence this equivariant infinite loop space structure on 
$\text{\bf BU}$ gives an equivariant spectrum $\gK\gR$ whose
associated  $RO(\bbz_2)$-graded cohomology theory is an enhancement of
$KR$-theory.  

In what follows, we show how to construct another $\bbz_2\ci_*$-functor
using constructions with algebraic cycles. The resulting equivariant
infinite loop space will then be used to provide {\it characteristic
classes} for the $RO(\bbz_2)$-graded $KR$-theory.
\smallskip

\Ex{6.8} {\bf (The algebraic cycle  functor)} \ \ Consider the 
functor defined by setting 
$T_Z(V) = \cz^v(\bbp(V_{\bbc}\oplus V_{\bbc}))$, 
the  topological group of  codimension-$v$ cycles in
$\bbp(V_{\bbc}\oplus V_{\bbc})$,
with $1_Z=1_G$. To a morphism $f:V \to W$ we associate
$$
T_Z(f)c \ = \
\bbp(f_{\bbc}(V_{\bbc})^{\perp}\oplus\{0\})\,\#\,
(f_{\bbc}\oplus f_{\bbc})_*(c), 
$$
and we define $\omega_Z$   by 
$$ 
\omega_Z(c,c') = \tau_*(c\# c').
$$
Using the same arguments as in [\LLM$_1$], it can be shown  that
$(T_Z, \omega_Z)$ is a $\bbz_2\ci_*$- functor with the $\bbz_2$-action
given by conjugation, and hence
$$
T_Z(\bbc^{\infty}) \ = \ \cz_{\bbc}^{\infty}
$$
is an equivariant  $\cl$-space.  In fact it is an equivariant
$E_{\infty}$-ring space which is additively $\bbz_2$-group complete and
therefore it is equivalent to the $0^{\text{th}}$  space of an
equivariant $E_{\infty}$-ring spectrum, which we denote by $\gZ_{\bbc}$.

The join operation $\omega_Z$ has various properties which yield
important results:
\roster
\item The join is multiplicative with respect to degree of cycles, in
other words $$\deg{\omega_Z(c,c')} = \deg{c}\cdot \deg{c'};$$
\item If $c$ is an averaged cycle then, for any fixed cycle $c'$, 
the join $\omega_Z(c,c')$ is also an averaged cycle. In other words,
the averaged cycles form an ``ideal'' within the fixed cycles.
\endroster

Now, let $\cz_{\bbc}^{\infty}(1)\subset \cz_{\bbc}^{\infty}$  be the
subspace consisting of the cycles of degree one. Since the join
operation $\omega_Z$ is equivariant and multiplicative on degrees, one
concludes that: 
\smallskip

\noindent{\bf a.}
$\cz_{\bbc}^{\infty}(1)$ is also an equivariant
$\bbz_2$-group-complete $\cl$-space, since it is connected and its
fixed point set is connected. It then follows that
$\cz_{\bbc}^{\infty}(1)$ carries a 
structure of an equivariant infinite loop space of its own, and hence
it is equivalent to the $0^{\text{th}}$  space of another
equivariant spectrum, which we denote by $\gM_{\bbr}$ to emphasize the
fact that we are equivariant delooping the {\it multiplicative units}
of $\cz_{\bbc}^{\infty}$.

\noindent{\bf b.}
Given $V \in \bbz_2\ci_*$, the ``forgetful map'' 
$$
\Phi_{V} : \cz^v_{\bbc}(\bbp(V_{\bbc}\oplus V_{\bbc}))\ \arr\ 
\cz^{2v}_{\bbr}(\bbp(V_{\bbc}\oplus \bar V_{\bbc}\oplus V_{\bbc}\oplus
\bar  V_{\bbc})) ,
$$
which sends c to $\omega_Z(c\oplus c)$, is not a group
homomorphism. Nevertheless, the preservation of degrees by the join
implies that the maps $\Phi_V$ define a map of  (non-equivariant)
$\ci_*$-functors between $\cz^*_{\bbc}$ and $\cz_{\bbr}^*$, which
preserves cycles of degree $1$. In particular, they induce a map of 
$\cl$-spaces $\Phi : \cz_{\bbc}^{\infty}(1) \to
\cz_{\bbr}^{\infty}(1).$ 
\smallskip

All of this discussion, in fact all the discussion in sections
4, 5, and 6 of  [\BLLMM] and section 3 of [\LLM$_1$],  
which include material on Chow monoid functors,  carries over directly
to our spaces of algebraic cycles. 

\Theorem{6.9}  {\sl  
\roster
\item The limiting topological group $\cz^{\infty}_{\bbc}$ is an
equivariant $E_{\infty}$-ring space which forms the $0$-level  space of
an equivariant $E_{\infty}$-ring spectrum $\gZ_{\bbc}$. The fixed point set
$\cz_{\bbr}^{\infty}$ is a (non-equivariant) $E_{\infty}$-ring space
which forms the $0$-level space of an $E_{\infty}$-ring spectrum
$\gZ_{\bbr}$.  The inclusion $\Gamma : \cz_{\bbr}^{\infty}
\hookrightarrow \cz_{\bbc}^{\infty}$ extends to a map of
(non-equivariant) ring  spectra $\Gamma : \gZ_{\bbr} \to \gZ_{\bbc}$.

\item The quotient group $\tilde{\cz}_{\bbr}^{\infty} \equdef
\cz_{\bbr}^{\infty} / (\cz_{\bbc}^{\infty} )^{av}$ is also an
$E_{\infty}$-ring space, and the quotient map 
$\rho : \cz_{\bbr}^{\infty} \to \tilde{\cz}_{\bbr}^{\infty}$
is a map of $E_{\infty}$-ring spaces. Hence, 
$\tilde{\cz}_{\bbr}^{\infty}$ is the $0^{\text{th}}$ space of an
$E_{\infty}$-ring spectrum $\tilde{\gZ}_{\bbr}$, and there is a
natural map of spectra $\rho : \gZ_{\bbr} \to \tilde{\gZ}_{\bbr}.$ 
Similarly, $\tilde{\cz}_{\bbr}^{\infty}(1)$ is an infinite loop space
under the operation induced by the join, which makes it into the
$0^{\text{th}}$ space of an spectrum $\tilde{\gM}_{\bbr}.$

\item $\cz^{\infty}_{\bbc}(1)$ carries an infinite loop space
structure which enhances the algebraic join, and makes it into the
$0$-level  space of an equivariant spectrum $\gZ_{\bbc}$. The fixed
point set $\cz_{\bbr}^{\infty}$ is a (non-equivariant)
$E_{\infty}$-ring space  which forms the $0$-level space of a
spectrum $\gM_{\bbr}$.  The inclusion   
$\Gamma : \cz_{\bbr}^{\infty}(1)  \hookrightarrow
\cz_{\bbc}^{\infty}(1)$    
extends to a map of  (non-equivariant) spectra 
$\gM_{\bbr} \to \gM_{\bbc}$. 

\item The canonical ``forgetful map'' 
$\Phi : \cz_{\bbc}^{\infty}(1) \to \cz_{\bbr}^{\infty}(1)$ 
induces a map of (non-equivariant) spectra 
$ \Phi : \gM_{\bbc} \to \gM_{\bbr}.$ 

\endroster } 
\medskip

An important feature of $\cz_{\bbc}^{\infty}(1)$ comes from the fact
that the   inclusion 
$$
G^v(\bbp(V_{\bbc}\oplus V_{\bbc})) \subset \cz^v(\bbp(V_{\bbc}\oplus
V_{\bbc} ))(1),
$$
as effective cycles of degree 1, is a natural transformation of
$\bbz_2\ci_*$-functors, and the resulting map 
$$
\text{\bf BU} \to \cz_{\bbc}^{\infty}(1)
$$
is an equivariant infinite loop space map. This fact,
together with the discussion above and [\BLLMM], gives the following result.
\medskip

\Theorem{6.10} {\sl 
\roster
\item The canonical equivariant inclusion 
$\text{\bf BU} \to \cz_{\bbc}^{\infty}(1)$  extends to a morphism 
$${ \bold{ \widetilde{c} } : \gK\gR \to \gM_{\bbc} }$$ 
of $\bbz_2$-equivariant spectra. 
Passing to fixed point sets gives maps of (non-equivariant) spectra 
$P : \gK\go \to \gM_{\bbr}$ from connective $KO$-theory to
$\gM_{\bbr}$, and $c : \gK\gu \to \gM_{\bbc}$ from connective
$K$-theory to $\gM_{\bbc}$. These maps fit into a commutative diagram
of spectra
$$
\CD
\gK\go @>{\gamma}>> \gK\gu \\
@V{{ \bold{\widetilde{c}}}}VV @VV{c}V \\
\gM_{\bbr} @>>{\Gamma}> \gM_{\bbc}
\endCD
$$
which extends the commutative diagram 
$$
\CD
\text{\bf BO} @>{\gamma}>> \text{\bf BU} \\
@V{P}VV @VV{c}V \\
\cz_{\bbr}^{\infty} @>>{\Gamma}>\cz_{\bbc}^{\infty},
\endCD
$$
where the map $c : BU \to \cz_{\bbc}^{\infty}$ classifies the {\rm
total Chern class.}
The composition $\gK\go \to \gM_{\bbr} \to \tilde{\gM}_{\bbr}$
is an extension to spectra level of the classifying map 
$\text{\bf BO}\to  \cz_{\bbr}^{\infty}(1) \to
\tilde{\cz}_{\bbr}^{\infty}(1)$ for the total Stiefel-Whitney class.
\endroster
}
\medskip

Analogous results for the groups of quaternionic cycles will be
established in the companion paper [\LLM$_3$].

A natural question now arises:  {\sl What is the equivariant 
cohomology theory classified by the $\bbz_2$-spectrum 
$\gZ_{\bbc}$}?  
In his thesis Pedro dos Santos has established the following beautiful
results. To state them we briefly recall some concepts from
equivariant homotopy theory (cf. [\May$_4$].)

Let $G$ be a finite group and $\underline M$ a Mackey functor for
$G$.  To each real representation $V$ of $G$ there is an
Eilenberg-MacLane space $K(\underline M , V)$ which classifies the 
ordinary equivariant cohomology group $H^V_G(\bullet ; \underline 
M)$ 
in dimension $V$ with coefficients in the Mackey functor 
$\underline M$.
These fit together to give an equivariant spectrum $\bbk(\underline
M,0)$ which classifies the full $RO_G$-graded equivariant 
cohomology
with coefficients in $\underline M$.

We now specialize to the group $G=\bbz_2$ and $\underline M= 
\underline {\bbz}$, the  Mackey functor constant at $\bbz$.  For each
$n$ we consider the fundamental respresentation $\bbr^{n,n} = 
\bbc^n$
of $\bbz_2$ given by complex conjugation.

\Theorem{6.11} {\bf ( dos Santos [dS])} {\sl There is a canonical
$\bbz_2$-equivariant homotopy equivalence
$$
\cz^{\infty}_{\bbc} \ \cong\ 
\prod_{n=0}^{\infty} K(\underline {\bbz}, \bbr^{n,n}).
\tag{6.1}
$$
This extends to an equivalence of $\bbz_2$-equivariant ring spectra
$$
\gZ_{\bbc} \ \cong\ 
\bbk(\underline {\bbz}, 0) \times \bbk(\underline {\bbz}, 
\bbr^{1,1})
\times \bbk(\underline {\bbz}, \bbr^{2,2}) \times \dots
$$
where $\bbk(\underline {\bbz}, 0)$ is the equivariant Eilenberg-
MacLane
spectrum and $\bbk(\underline {\bbz}, \bbr^{n,n})$ is the connective 
 equivariant spectrum with 
$\Omega^{\bbr^{n,n}}\bbk(\underline {\bbz}, \bbr^{n,n}) \cong
\bbk(\underline {\bbz}, 0)$, and where the ring structure is given by
the equivariant cup product pairing. }

\bigskip

For a $\bbz_2$-space $X$ we denote by 
$H^*_{\bbz_2}(X;\underline{\bbz})$ the full $RO_{\bbz_2}$-graded
equivariant cohomology ring of $X$ with coefficients in the Mackey
functor $\underline {\bbz}$. We  abbreviate 
$H_{\bbz_2}^{\bbr^{n,n}}(X;\underline{\bbz}) \equiv
H_{\bbz_2}^{n,n}(X;\underline{\bbz})$

\medskip\ 

\Theorem{6.12} {\bf (Dugger   and  dos Santos  )} {\sl  There is a
canonical ring homomorphism 
$$
H^*_{\bbz_2}({\bold {BU}}; \underline {\bbz}) \cong 
R[{\tc}_1, {\tc}_2, {\tc}_3,   \dots ]
$$
where
$$
{\tc}_n \in H_{\bbz_2}^{n,n}({\bold {BU}}; \underline
{\bbz}) 
$$
for each $n$ and  
$
R = H^*_{\bbz_2}({\bold {pt}}; \underline {\bbz})
$
 is the coefficient ring.

Furthermore, let $\widetilde{\iota}_{n,n}$ denote the fundamental
class of  $ K(\underline {\bbz}, \bbr^{n,n})$.  Then with respect
to the splitting (6.1) the natural $\bbz_2$-map 
$$
P :{\bold{BU}} \ \arr\ \cz^{\infty}_{\bbc} 
$$
satisfies
$$
P^*\left(\widetilde{\iota}_{n,n}\right) = {\tc}_n
$$

}

\Note { } The first assertion of Theorem 6.12 is due to Dan
Dugger [D] and the second to dos Santos [dS].

\medskip

Theorem 6.12 shows that the inclusion map 
$\tc :  {\bold{BU}} \to \cz^{\infty}_{\bbc}(1)$
naturally classifies the {\sl total Chern class in (full
$RO_{\bbz_2}$-graded) equivariant cohomology}.  Thus for Real 
spaces
$X$, $\tc$ determines a natural transformation
$$
\tc : KR(X) \ \arr\ 
\bigoplus_{n\geq 0} H_{\bbz_2}^{n,n}(X;\, \underline \bbz)
$$
and the property that $\tc (V\oplus V') = 
\tc (V) \# \tc (V')$, together with Theorem 6.12,
shows that
$$
\tc (E\oplus E') \ =\ \tc (E) \cup \tc (E')
$$
for all $E,E' \in KR(X)$.  Theorem 6.10 shows that the equivariant
infinite loop structure on $\cz_{\bbc}^{\infty}(1)$ corresponding to
the spectrum $\gM_{\bbc}$ makes this total Chern class map 
$\tc $ {\sl an equivariant infinite loop map}.  This is the 
full $\bbz_2$-equivariant version of   Segal's conjecture settled in
[\BLLMM].

In fact under the forgetful functor the class $\tc$
becomes the ordinary total Chern class $c$, and the map
$ {\bold \tc}: \gK\gR \to \gM_{\bbc}  $ of equivariant
spectra becomes the map
$ {\bold  c}: \gK\gu \to \gM_{\bbc}  $
of non-equivariant spectra studied in [\BLLMM].

More interesting perhaps is the restriction of $\tc$ to the 
fixed-point set ${\bold{BO} \subset \bold{BU}}$.  This gives a
characteristic class for real bundles which is a mixture of 
$\bbz$ and $\bbz_2$ classes and satisfies Whitney duality.  We 
examine this next.

\subheading{\S 7.  A new total characteristic class  }
The mappings 
$$
P : G^q_{\bbr}(\bbp(V)) \ \arr\ \cz^q_{\bbr}(\bbp(V)) 
$$
which stabilize to 
$$
P: \BO_q \ \arr\ \cz^q_{\bbr}
\qquad\text{and}\qquad
P : \text{\bf BO}\arr\ \cz^{\infty}_{\bbr}
$$
represent a   ``total'' characteristic class which is in complete 
analogy
with the {\sl total Chern class}: 
$$
c : G^q_{\bbc}(\bbp(V)) \ \arr\ \cz^q_{\bbc}(\bbp(V)) 
$$
$$
c: \BU_q \ \arr\ \cz^q_{\bbr}
\qquad\text{and}\qquad
c : \text{\bf BU}\arr\ \cz^{\infty}_{\bbr}
$$
and the {\sl total Stiefel-Whitney class}: $$
w : G^q_{\bbr}(\bbp(V)) \ \arr\ \widetilde{\cz}^q_{\bbr}(\bbp(V)) 
$$ 
$$
w: \BO_q \ \arr\ \widetilde{\cz}^q_{\bbr}
\qquad\text{and}\qquad
w : \text{\bf BO}\arr\ \widetilde{\cz}^{\infty}_{\bbr}.
$$
This new class has the property that for real vector bundles $E$ and 
$F$
over a space $X$,
$$
P(E\oplus F) \ = \ P(E)\cup P(F).
\tag{7.1}
$$
As a map from $\text{\bf BO}$ to ${\cz}^{\infty}_{\bbr}$ it is an 
infinite loop
map.  It fits into a pattern of infinite loop diagrams as we saw in \S 
6.
This compelling picture makes the  further study of $P$ intriguing.

\medskip

Recall that we have a commutative diagram:
$$\CD
\ @. \  @. \ @. \ @.   \cz^{\infty}_{\bbc}  \\
\ @.\   @. \nearrow @. \ @. \      \\
\text{\bf BO} @>{P}>> {\cz}^{\infty}_{\bbr} @. \ @. \ @. \  \\
\ @. \ @. \searrow @. \ @. \  \\
\ @. \  @. \ @. \ @.  \widetilde{\cz}^{\infty}_{\bbr}  \\
\endCD$$
where the upper map represents the total Chern class of the 
complexification
(essentially the total Pontrjagin class) and the lower map is the total
Steifel-Whitney class.  However, $P$ contains much more 
information.
From the splitting in Theorem 3.3 the map $P$ is seen to represent a 
certain sum
of integral and mod 2 cohomology classes.  Thus $P$ is a particular 
arrangement
of Pontryagin and Stiefel-Whitney classes.  Exactly which 
arrangement 
has been determined by dos Santos.

\Theorem {7.1}  ([\dS])  \ {\sl  For $k<n$, one has
$$
P^*(\iota_{n,k}) = Sq^k w_n
$$
where $Sq^k$ denotes the $k^{\text{th}}$ Steenrod operation and 
$w_n$ denotes the $n^{\text{th}}$ Stiefel-Whitney class.

When $k=n$ and $n$ is even, $P^*(\iota_{n,n}) $ is the $n^{\text{th}}$
Pontrjagin class $p_n$.
}

There is a second construction that one
can   associate to Real bundles.

\medskip\noindent
{\bf Construction 7.2} \  To a Real map $f:X \to G^{q}(\bbp(\bbc^n))$
classifying a Real bundle $E_f\to X$, we associate the  mapping
$$
\tf : X/\bbz_2 \ \arr\ \cz^{q}(\bbp(\bbc^N))^{av}  
$$
defined by $\tf([x]) = f(x) + f(\rho x) = f(x) + \rho f(x)$.  Let 
$\iota _{n,k}$ denote the fundamental class of the factor $K(I_{n,k}, 
n+k)$
in the canonical splitting of $\cz^q_{\bbr}$ given in Theorem 3.3. 
Then
$$
\tf^*(\iota _{n,k}) \in H^*(X/\bbz_2)
$$
is an invariant of the Real bundle $E_f$.

\Ex{7.3}  \ \ Consider the commutative diagram of Real spaces
$$\CD
\bbp^n @>>> G^n(\bbp(\bbc^N)) \\
@VVV   @VVV  \\
\cz_0(\bbp^n)  @>>> \cz^n(\bbp(\bbc^N))
\endCD
$$
where the left vertical map is the standard inclusion, the horizontal 
maps are
complex suspension, and the right vertical map is our $\bbz_2$-
characteristic
map.
Let $\bbp^{n,k}\subset \bbp^n$ be the Real spaces which are 
introduced in the  proof of
Theorem 9.1.  Composing with suspension gives a Real mapping 
$
\bbp^{n,k} \arr G^n(\bbp(\bbc^N)),
$
for which 7.1 yields a map 
$
\tf: \bbp^{n,k}/\bbz_2 \arr \cz^{q}(\bbp(\bbc^N))^{av}.
$
It follows from the discussion (9.11) ff. that the corresponding 
classes
$\tf^*(\iota _{n,k})$ are non-trivial for $0<k \leq n$.

\Note{7.4}  Let $\widetilde{P}(E)$ denote the total class
of a Real bundle $E$ constructed in 7.2.  Then we have 
$$
\widetilde{P} (E\oplus E') \ = \ \widetilde{P(E)\# P(E')}
$$
This class satisfies the addition relation:
$$
\widetilde{P} (E\oplus E') + \widetilde{P} (E\oplus \overline {E}')
\ =\ \widetilde{P} (E)\#\widetilde{P} (E')
$$

 \heading \S 8. Proof of Theorems 3.3 and 3.4  
 \endheading

\define\N{q}
\define\n{n}
\define\k{k}
\define\M{N}

Consider the Real vector space $(\bbc^{n+1}, \tau)$ where 
$\tau$ is complex
conjugation, and denote the corresponding projective space by 
$\pc^n$ and its real
form (the $\tau$-fixed-point set) by $\pr^n$.
We will make  thorough use of the identification of
$\pc^n$  with the $n$-fold symmetric 
product $ \SP{n}{\pc^1}$ of $\pc^1$, and the fact that the complex
conjugation  involution on $\pc^n$ is induced by the complex
conjugation $\tau$ on $\pc^1$ under this identification. 

\Prop{8.1}\ {\sl
Fix integers $\N \leq \M$.
Then
\roster
\item "{\bf a.}" Iterations of the algebraic suspension give
 canonical homotopy  equivalences:
$$ 
\alignat {1} 
\crl{\N}{\pc^\M}  & \cong \crl{\N}{\pc^{\N}}  = \cycd{0,\bbr}{\pc^{\N}}, \\
\cyav{\N}{\pc^\M} & \cong \cyav{\N}{\pc^{\N}} = \cycd{0}{\pc^{\N}}^{av}, \\ 
\crd{\N}{\pc^\M}  & \cong \crd{\N}{\pc^{\N}}  = \crdd{0}{\pc^{\N}}.
\endalignat
$$

\item "{\bf b.}" For any compact space $X$ with an action of $\bbz_2$
there is a natural degree-preserving topological isomorphism:
$\cycd{0}{X}^{av} \cong \cycd{0}{X/\bbz_2},$
where the latter denotes the free abelian group on the space
$X/\bbz_2$.
In particular, there is a natural topological isomorphism
$\cyav{\N}{\pc^\N} \cong \cyc{\N}{\pc^\N/\bbz_2}=
\cycd{0}{\pc^\N/\bbz_2}$.

\item "{\bf c.}" The  short exact sequence
$$ 0\rightarrow \cyav{\N}{\pc^\N} \rightarrow
   \crl{\N}{\pc^\N}
   \rightarrow \crd{\N}{\pc^\N} \rightarrow 0 .
$$
is a principal  fibration.
\endroster
}
\pf
Part {\bf a} follows from repeated applications of Theorem 3.2. 
Now, consider the topological homomorphism 
$\psi : \cycd{0}{X} \rightarrow \cycd{0}{X}^{av} \subset \cycd{0}{X}$ 
defined by 
$\psi(\sigma) = \sigma + \tau*\sigma$, 
where  $\tau*\sigma$ denotes the action of the generator of $\bbz_2$
on $\sigma$. Since $X$ is compact, it follows from the description of
the topology of $\cycd{0}{X}$ that $\psi$ is a closed map, which 
clearly surjects onto $\cycd{0}{X}^{av}$.
The composition $X\rightarrow \cycd{0}{X} \rightarrow
\cycd{0}{X}^{av}$ clearly factors through the projection 
$\pi: X\rightarrow X/\bbz_2\ $, and hence the universal property of the
functor $\cycd{0}{-}$ gives a topological homomorphism
$\Psi : \cycd{0}{X} \rightarrow \cycd{0}{X}^{av}$ 
such that 
$\Psi \circ\pi_* = \psi$, where 
$\pi_* : \cycd{0}{X} \rightarrow \cycd{0}{X/\bbz_2}$ is the projection
induced by $\pi$. It is a routine verification to see that $\Psi $ is
injective, and hence a closed continuous bijection. This proves part {\bf b}.

In order to prove part {\bf c}, consider the monoid
$C \equdef \amalg_{d\geq 0}\  SP_d(\pc^q)^{fix}$ and its closed submonoid 
$C' \equdef \amalg_{d'\geq 0}\ SP_{2d'}(\pc^q)^{av}$. Note that the 
na\"{\i}ve group completions of $C$ and $C'$ are 
$\crl{q}{\pc^q}$ and $\cyav{q}{\pc^q}$, respectively. Since complex
conjugation induces a real analytic map on all 
products $SP_d(\pc^q) \times SP_{d'}(\pc^q)$, preserving filtrations
by degrees, one can provide equivariant triangulations to all such
products, making $(C,C')$ into a triangulated pair. It follows
that $(C,C')$  satisfies the hypothesis of 
[\Li$_2$, Theorem 5.2], which then implies the desired result.
\qed

\Cor{8.2}\ {\sl For any space $X$ with a $\bbz_2$ action, there 
is a  natural isomorphism
$$
\pi_j(\cyavd{0}{X}) \cong H_j(X / \bbz_2; \;
\bbz).
\tag{8.1}
$$
Therefore, one has a canonical isomorphism
$$
\pi_j(\cyav{\N}{\pc^\N}) \cong H_j(\pc^\N / \bbz_2; \;
\bbz),
\tag{8.2}
$$
and the homotopy type
of $\cyav{\N}{\bbp^\N_{\bbc}}$ is completely determined by the
singular homology of the quotient space $X^\N = \pc^\N / \bbz_2$.}
\medskip
\pf
It follows from Proposition 8.1(c) that
$\pi_j(\cycd{0}{X}^{av}) \cong \pi_j(\cycd{0}{X/\bbz_2})$, whereas the
latter group is naturally isomorphic to $H_j(X/\bbz_2; \bbz)$ by the
Dold-Thom theorem. This proves (8.1). The last assertion follows from
the first one, and the fact that $\cycd{q}{\pc^q}^{av}$ is a product of
Eilenberg-MacLane spaces, since it is an abelian topological group.
\qed
\medskip

The first fundamental ingredient, for our splitting results, is
given by the following equivariant version of certain constructions
in [\FL].  

Let
$x_0 \in \pc^1$ be a basepoint which is fixed under the complex
conjugation $\tau$, i.e. a real point. Then, for $q\leq n$, the
canonical inclusion 
$$
\aligned
i_{\n,\N} : \pc^\n = \SP{\n}{\pc^1}   & \longrightarrow \pc^\N =
\SP{\N}{\pc^1}  \\
\sigma & \longmapsto \sigma + (\N-\n)x_0
\endaligned
 \tag{8.3}
$$
is $\bbz_2$-equivariant.
Now, given $\n\leq \N$, define a map
$r_{\N,\n} : \pc^\N \to \SP{\binom{\N}{\n}}{\pc^\n}$   by sending
$ a_1+\cdots+a_\N \in \SP{\N}{\pc^1}\equiv \pc^\N$  to
$ \sum_{ |I|= \n } \{ a_{i_1}+\cdots+a_{i_\n} \} \in
\SP{\binom{\N}{\n}}{\SP{\n}{\pc^1}} 
\equiv\SP{\binom{\N}{\n}}{\pc^\n}  .$
The $\bbz_2$ action on $\pc^\n$, induced by conjugation, extends
``linearly'' to a $\bbz_2$ action on $\SP{\binom{\N}{\n}}{\pc^\n}$ in 
such
a way that the map $r_{\N,\n}$ becomes $\bbz_2$-equivariant.
This map, in turn, induces an equivariant  homomorphism
$$
\br_{\N,\n} : \cycd{0}{\pc^\N} \to
\cycd{0}{\SP{\binom{\N}{\n}}{\pc^\n}}.
\tag{8.4}
$$

One has an evident ``trace map'' (see [Proposition 7.1, \FL$_1$])
$$
tr : \cycd{0}{\SP{\binom{\N}{\n}}{\pc^\n}} \to
\cycd{0}{\pc^\n}
\tag{8.5}
$$
defined as the extension to the free abelian group  
$ \cycd{0}{\SP{\binom{\N}{\n}}{\pc^\n}}$
of the natural inclusion
$ \SP{\binom{\N}{\n}}{\pc^\n} \hookrightarrow \cycd{0}{\pc^\n}$.
This is easily seen to be $\bbz_2$-equivariant and
can be used to define an equivariant homomorphism
$$
\rho_{\N,\n} : \cycd{0}{\pc^\N} \to \cycd{0}{\pc^\n}
\tag{8.6}
$$
as the composition
$\cycd{0}{\pc^\N} @>{\br_{\N,\n}}>>
\cycd{0}{\SP{\binom{\N}{\n}}{\pc^\n}} @>{tr}>>
\cycd{0}{\pc^\n}$. Note that $\rho_{\n,\n}$ is 
the identity map.

Now, let
$\psi_\n : \cycd{0}{\pc^\n} \to
\cycd{0}{\pc^\n}/\cycd{0}{\pc^{\n-1}}$
denote the quotient map and observe that this is an equivariant map,
for the induced $\bbz_2$ action on the quotient group.
Let
$$
q_{\N,\n} : \cycd{0}{\pc^\N} \to
\cycd{0}{\pc^\n}/\cycd{0}{\pc^{\n-1}}
\tag{8.7}
$$
denote the composition
$\psi_\n \circ \rho_{\N,\n}$. We adopt the convention that $\pc^{-1} 
= \emptyset,$
and $\cycd{0}{\emptyset} = \{ 0 \}.$

Since $\psi_\n$, $i_{\n,\N}$ and $q_{\N,\n}$ are all
equivariant maps, they induce homomorphisms on the level of fixed 
point
set and averaged  subgroups, whose properties are described in the
following result.

\Prop {8.3}\ {\sl
Let $\{ ( M_\N, Q_\n, q_{\N,\n}, i_{\n,\N},  ) \ | \ q_{\N,\n} : M_\N 
\to Q_\n,
\ \ i_{\n,\N} : M_\n \hookrightarrow M_\N,\ \  0\leq \n\leq \N \}$
denote  one of the following collections of spaces and maps:
\roster
\item  "{\bf a.}"\ $M_\N = \cycd{0}{\pc^\N}$, $Q_\n = 
\cycd{0}{\pc^\n}/\cycd{0}{\pc^{\n-1}}$,
with $q_{\N,\n}$ and $i_{\n,\N}$ given by (8.7) and (8.3), 
respectively.
\item   "{\bf b.}"\ $M_\N = \crld{0}{{\pc^\N}}$, $Q_\n = \{
\cycd{0}{\pc^\n}/\cycd{0}{\pc^{\n-1}} \}^{fix}$,
with $q_{\N,\n}$ and $i_{\n,\N}$ induced  by (8.7) and (8.3) on the
respective fixed point sets.
\item   "{\bf c.}"\ $M_\N = \cyavd{0}{\pc^\N}$, $Q_\n = \{ \cycd{0}{\pc^\n}/
\cycd{0}{\pc^{\n-1}} \}^{av}$,
with $q_{\N,\n}$ and $i_{\n,\N}$ induced by (8.7) and (8.3) on the
respective subgroups of averaged cycles.
\endroster
Then the following assertions hold:
\roster
\item The sequence $M_{\n-1} @>{i_{\n-1,\n}}>> M_\n 
@>{q_{\n,\n}}>> Q_\n$ of
topological group homomorphisms is a principal fibration, for all 
$\n$.
\item The following diagram commutes:
$$
\CD
M_\n @>{i_{\n,\N}}>> M_\N \\
@V{q_{\n,\n}}VV @VV{q_{\N,\n}}V \\
Q_\n @>>{id}> Q_\n.
\endCD
$$
\endroster
}
\pf
Since complex conjugation is a real analytic map, the pair
$(\pc^n,\pc^{n-1})$ becomes a $\bbz_2$-simplicial pair after a
suitable equivariant triangulation. Therefore, assertion (1) for the
collections {\bf a} and {\bf b} above follows from 
[\Li$_4$, Theorem 2.7]. Using Proposition 8.1(b), the exact sequence 
$\cycd{0}{\pc^{n-1}}^{av}\rightarrow \cycd{0}{\pc^{n}}^{av}
\rightarrow \cycd{0}{\pc^{n}}^{av}/\cycd{0}{\pc^{n-1}}^{av}$ 
becomes
$$\cycd{0}{\pc^{n-1}/\bbz_2}\rightarrow \cycd{0}{\pc^{n}/\bbz_2}
\rightarrow \cycd{0}{\pc^{n}/\bbz_2}/\cycd{0}{\pc^{n-1}/\bbz_2}.
\tag{8.8}
$$
It is well known that, since $(\pc^{n}/\bbz_2, \pc^{n-1}/\bbz_2)$
is a CW-pair, then (8.8) is a principal fibration; cf. [\DT] or [\Li$_2$].
\define\mby {\bold y}
\define\mbx {\bold x}

To prove assertion (2), first observe that all maps in the diagram for
collection {\bf a} are $\bz_2$-equivariant. Therefore, assertion (2) for
collections {\bf b} and {\bf c} follows from its validity for the collection
{\bf a}. 

Given $\mbx = x_1+\cdots + x_n \in \pc^n \equiv SP^n(\pc^1)$,
denote $\mby = i_{n,q}(\mbx)$, and observe that $\mby = y_1+\cdots y_q$, 
where $y_i = x_i$ for $i=1,\ldots,n$, and $y_i = x_0$ for
$i=n+1,\ldots,q$. It follows from its definition that 
$r_{q,n}(\mby) = \sum_{|I|=n} \{ y_{i_i}+\cdots +y_{i_n} \}$
can be written as $r_{q,n}(\mby) = \{ x_1+\ldots x_n\} + 
\sum_{|J|=n-1} \{ x_0 + y_{j_1} +\cdots + y_{j_{n-1}} \}$.
Therefore we obtain
$$
\alignat { 1 }
q_{q,n} \circ i_{n,q} (\mbx) & = 
\psi_n \circ \rho_{q,n}\circ i_{n,q} (\mbx) = 
\psi_n \circ tr \circ r_{q,n}(\mby) \\
& = \psi_n \circ tr \left( \{ x_1+\ldots + x_n\} + 
\sum_{|J|=n-1} \{ x_0 + y_{j_1} +\cdots + y_{j_{n-1}} \} \right) \\
& = \psi_n \left( \{ x_1+\ldots x_n\} + 
\sum_{|J|=n-1} \{ x_0 + y_{j_1} +\cdots + y_{j_{n-1}} \} \right) \\
& = \psi_n \left( \{ x_1+\ldots x_n\}  \right).
\endalignat
$$
The latter equality comes from the fact that
$\sum_{|J|=n-1} \{ x_0 + y_{j_1} +\cdots + y_{j_{n-1}} \} $
is supported in $\pc^{n-1}\subset \pc^{n}$ and hence lies in the
kernel of $\psi_n$. We have then shown that 
$q_{q,n}\circ i_{n,q} (\mbx) = \psi_n (\mbx) = \psi_n \circ \rho_{n,n}
(\mbx) = q_{n,n}(\mbx)$, since $\rho_{n,n}$ is the identity and 
$q_{n,n}\equdef \psi_n\circ \rho_{n,n}$. Since $M_n = \cycd{0}{\pc^n}$
is freely generated by the points $\mbx \in \pc^n$, the result follows.
\qed

\Cor {8.4}\ {\sl
If  $ M_\N, Q_\n, q_{\N,\n}$ and $i_{\N,\n}$ are as above, then the 
group
homomorphism
$$
M_\N \rightarrow Q_0 \times \cdots \times Q_\N
$$
given by the product $q_{\N,0}\times \cdots \times q_{\N,\N}$
is a homotopy equivalence.
}
\pf
The proposition above guarantees that all three collections considered
satisfy the hypothesis of [\FL$_1$, Proposition 2.13], from which the
corollary follows.
\qed
\smallskip

Recall that, given a $CW$-pair $(X,Y)$ with $X$ connected and a
basepoint $x_0 \in Y$, then the quotient group
$\cycd{0}{X}/\cycd{0}{Y}$ is naturally identified with 
$\cycd{0}{X/Y}_o$,
the connected component of $0$ in $\cycd{0}{X/Y}$. Furthermore, 
this connected component is the subgroup consisting of cycles of degree
$0$; cf. [\Li$_4$, \S 2]. For uniformity of notation, we write
$\cycd{0}{X / \emptyset}_o = \cycd{0}{X}$.

The canonical  homotopy equivalences
$$
\cyav{\N}{\pc^\M} \simeq \cyav{\N}{\pc^\N} \simeq 
\prod_{\n=0}^\N \{
\cycd{0}{\pc^\n}/\cycd{0}{\pc^{\n-1}} \}^{av},
\tag{8.9}
$$
respectively described in Propositions 8.1(a) and Corollary 8.4,
together with  the previous paragraph, give a homotopy equivalence
$$
\cyav{\N}{\pc^\M} \simeq \prod_{\n=0}^\N\ \{ \cycd{0}{\pc^\n / 
\pc^{\n-1} }_o \}^{av},
\tag{8.10}
$$
and Proposition 8.1(c) then gives
$$
\cyav{\N}{\pc^\M} \simeq \prod_{\n=0}^\N \cycd{0}{ \{ \pc^\n / 
\pc^{\n-1} \}/ \bbz_2 }
_o.
\tag{8.11}
$$
On the other hand, it is easy to see that one has homeomorphisms
$$
\{ \pc^{\n}/ \pc^{\n-1}\}/ \bbz_2 \cong S^{2\n}/\bbz_2 \cong S^\n 
\# \pr^{\n-1},
\tag{8.12}
$$
where $S^{2\n}$ is the one point compactification of $\bbc^\n$, and
$\#$ denotes the real join of topological spaces. The
$\bbz_2$ action  on $S^{2\n}$ here is the extension of the usual
complex conjugation on $\bbc^\n$.

Given $\n\geq 1$ one has a natural isomorphism
$\tilde{H}_\k(S^\n \# \pr^{\n-1}; \bbz)\cong
\tilde{H}_{\k-\n-1}(\pr^{\n-1};\bbz)$, and hence
the spaces $\{ \pc^{\n}/ \pc^{\n-1}\}/ \bbz_2$
satisfy the hypothesis of Theorem A.5. 
As a consequence one obtains a canonical splitting 
$$
\cycd{0}{\{ \pc^{\n}/ \pc^{\n-1}\}/ \bbz_2}_o \cong
\prod_{\k=0}^{2\n} K(\tilde{H}_\k(S^\n\# \pr^{\n-1};\bbz),\ \k) 
\cong
\prod_{\k=n+1}^{2\n} K(\tilde{H}_{\k-\n-1}( \pr^{\n-1};\bbz),\ \k).
\tag{8.13}
$$

The following result synthesizes the above conclusions.

\Theorem{8.5}\ {\sl
The group of averaged cycles of codimension $\N$ and degree $0$ in
$\pc^\M$ is connected and has a canonical splitting into products of
Eilenberg-MacLane spaces
$$
\cyav{\N}{\pc^{\M}} \cong \prod_{\n=0}^\N \prod_{\k=1}^{\n}
K\left( {H}_{\k-1}(\pr^{\n-1};\bbz),\ \k+\n \right).
$$
The homotopy groups
$\pi_*(\cyav{\N}{\pc^{\M}})$
have the  structure of a bigraded abelian group
$$
 \pi_*(\cyav{\N}{\pc^{\M}})\cong \bigoplus_{\n,\k\geq 0} 
\Iav{\n}{\k},
\quad 
\text{where}  
\quad 
\Iav{\n}{\k} \equiv \pi_{\n+\k}\left\{\cz_0(\pc^{\n}/\pc^{\n-
1})_o\right\}^{av} 
\cong \tilde{H}_{\k-1}(\pr^{\n-1};\bbz),
$$
 for $\n+\k>0$, and $I_{0,0} = 2\bbz$.
In other words,
$$\Iav{\n}{\k} = \cases
            2\bbz &, \text{ if } \n=\k=0 ;\\
            0 &, \text{ if $\k$ is odd,   or $\k>\n$, or\   } \n> \N; \\
            \bbz &, \text{ if $\k=\n\leq \N$ and $\k\geq 2$ is even; } \\
            \bbz_2 &, \text{ if $\k< \n \leq \N$  and  $\k\geq 2$ is 
even}.
             \endcases
$$
}
\smallskip
\pf
The splitting map is just the combination of the canonical splittings
(8.9) through (8.13). The bigraded group structure is
self-explanatory, except for the identity $\Iav{0}{0} =2\bbz$. This
stems from the fact that the restriction of the degree homomorphism
$\deg : \cycd{0}{\pc^{\N}} \to \bbz$ to 
$\cycd{0}{\pc^\N}^{av}$ induces an isomorphism between 
$\pi_0(\cycd{0}{\pc^\N}^{av})$ and $2\bbz \subset \bbz$. Indeed, given
$\sigma \in \cycd{0}{\pc^\N}^{av}$, one can write $\sigma = u +
\tau*u$. If $x_o \in \pc^\N$ is the basepoint and $d = \deg{u}$, then
one can find a path $\gamma $ in $\cycd{0}{\pc^\N}$ between
$u$ and $d\cdot x_o$. The path $\gamma + \tau*\gamma$ lies in
$\cycd{0}{\pc^\N}^{av}$ and connects $\sigma$ to $2d\cdot x_o$. 
\qed
\medskip

\Remark {8.6}\ Note that Theorem 3.4 is simply a reformulation of
Theorem 8.5.
\medskip

One of the pleasant consequences of Proposition 8.1, and the arguments
in its proof, is the fact that one has a commutative diagram
$$
\CD
\cyavd{0}{\pc^\N} @>{sp^{av}}>>
\prod_{\n=0}^\N\ \{ \cycd{0}{\pc^\n / \pc^{\n-1} }_o\}^{av} \\
@V{i}VV @VV{\prod_\n i_\n }V \\
\crld{0}{\pc^\N} @>{sp^{fix}}>> \prod_{\n=0}^\N\ \{ \cycd{0}{\pc^\n / 
\pc^{\n-1}
}_o\}^{fix} \\
@V{q}VV @VV{\prod_\n q_\n}V \\
\crdd{0}{\pc^\N} @>{sp^{red}}>> \prod_{\n=0}^\N\ \{ \cycd{0}{\pc^\n 
/ \pc^{\n-1}
}_o\}^{fix} /  \{ \cycd{0}{\pc^\n / \pc^{\n-1} }_o\}^{av} .
\endCD
\tag{8.14}
$$
whose columns are fibrations. It then follows from the  five lemma
that the bottom horizontal arrow is also a homotopy equivalence. 

Notice that one has canonical topological isomorphisms
$$
\aligned
\{ \cycd{0}{\pc^\n / \pc^{\n-1}}_o\}^{fix}/
\{ \cycd{0}{\pc^\n /\pc^{\n-1}}_o \}^{av} &\cong
\{ \cycd{0}{S^{2\n}}_o \}^{fix} /  \{ \cycd{0}{S^{2\n} }_o \}^{av}  \\
& \cong \cycd{0}{ \{ S^{2\n} \}^{fix} }_o\otimes \bbz_2 \cong
\cycd{0}{ S^{\n} }_o \otimes \bbz_2 ,
\endaligned
\tag{8.15}
$$
for $\n\geq 1.$
This provides the canonical splitting of Theorem 3.5 [\Lam]:
$$
\crdd{0}{\pc^\N} \cong \prod_{\n=0}^\N \cycd{0}{S^{\n}}_o\otimes 
\bbz_2
\equiv \prod_{\n=0}^\N K(\bbz_2,\n).
\tag{8.16}
$$

It follows that each factor in the right hand column of (8.14)
yields a long exact sequence of homotopy groups
$$
\aligned
\cdots & \to
\pi_{\k+1}(\cycd{0}{S^\n}_o\otimes \bbz_2 ) \to
\pi_{\k}(\cycd{0}{S^\n\# \pr^{\n-1}}_o) \to
\pi_{\k}(\{ \cycd{0}{S^{2\n}}_o\}^{fix} ) \to  \\  &\to
\pi_{\k}(\cycd{0}{S^\n}_o\otimes \bbz_2 )  \to
\pi_{\k-1}(\cycd{0}{S^\n\# \pr^{\n-1}}_o) \to
\cdots.
\endaligned
$$
Therefore, for $\k \neq \n-1,\ \n$ the inclusion
$$
i_\n : \cycd{0}{\pc^\n/\pc^{\n-1}}_o^{av}= \cycd{0}{S^\n\# \pr^{\n-
1}}
\to
\cycd{0}{\pc^\n/\pc^{\n-1}}_o^{fix} =\cycd{0}{S^{2\n}}_o^{fix}
$$
induces an isomorphism in homotopy groups $\pi_\k$. It remains to 
examine the     
exact sequence 
$$
\aligned
0 & \to
\pi_{\n}(\cycd{0}{S^\n\# \pr^{\n-1}}_o) \to
\pi_{\n}(\{ \cycd{0}{S^{2\n}}_o\}^{fix} ) @>{q_{\n*}}>>
\pi_{\n}(\cycd{0}{S^\n}_o\otimes \bbz_2 )  \to \\
& \to \pi_{\n-1}(\cycd{0}{S^\n\# \pr^{\n-1}}_o) \to
\pi_{\n-1}(\{ \cycd{0}{S^{2\n}}_o\}^{fix} ) \to 0.
\endaligned
$$
Since
$\pi_{\n}(\cycd{0}{S^\n\# \pr^{\n-1}}_o) \cong
 \pi_{\n-1}(\cycd{0}{S^\n\# \pr^{\n-1}}_o) \cong 0,$
one concludes that
the quotient
$q_\n : \{ \cycd{0}{S^{2\n}}_o\}^{fix} \to
\cycd{0}{S^\n}_o\otimes \bbz_2
$
induces an isomorphism on the $\n$-th homotopy group:
$$
 \pi_{\n}(\{ \cycd{0}{S^{2\n}}_o\}^{fix} )   \cong  
\pi_{\n}(\cycd{0}{S^\n}_o\otimes \bbz_2 ) \cong \bbz_2
\quad\text{and}\quad
\pi_{\n-1}(\{ \cycd{0}{S^{2\n}}_o\}^{fix} ) = 0.
\tag{8.17}
$$
We have thus completed the computation of the homotopy type of 
the space
of real cycles on projective spaces.

\Theorem{8.7}\ {\sl
The homotopy groups of $\crl{\N}{\pc^{\M}}$
have the  structure of a bigraded abelian group
$$
\pi_*(\crl{\N}{\pc^{\M}})\cong 
\bigoplus_{\n,\k\geq 0} \I{\n}{\k},\ \ \text{where}\ \ 
\I{\n}{\k} \equiv 
\pi_{\n+\k}\left\{\cz_0(\pc^{\n}/\pc^{\n-1})_o\right\}^{fix}
$$
and
$$\I{\n}{\k} = 
            \cases
            0 &,  \text{ if $\k$ is odd, or  } \k>\n,\text{  or } \n> \N;\\
            \bbz &, \text{ if $\k=\n$ and $\k$ is even; } \\
            \bbz_2 &, \text{ if $\k < \n \leq \N$ and  $\k\geq 0$ is 
even}.
             \endcases
$$
}
\pf
The computation of the homotopy groups comes from Theorem 8.5, 
(8.13) and (8.14). However, we need to explain once again the
identification $\I{0}{0} =\bbz$, which is once again given by the
restriction of the degree map to $\cycd{0,\bbr}{\pc^\N}$. Indeed, any
cycle  $ \sigma \in\cycd{0,\bbr}{\pc^\N} $ can be written as a sum 
$\sigma = u+v$, where 
$u\in \cycd{0}{\pr^\N}\subset \cycd{0,\bbr}{\pc^\N}$ has degree $d$
and  $u\in \cycd{0}{\pc^\N}^{av}\subset \cycd{0,\bbr}{\pc^\N}$ has
degree $2e$. Since $\pr^\N$ is connected, and  $x_o\in \pr^\N \subset \pc^N$  
one finds paths $\gamma$ in $\cycd{0}{\pr^\N}$ and $\tilde{\gamma}$ in
$\cycd{0}{\pc^\N}^{av}$ connecting $u$ to $d\cdot x_o$ and $v$ to
$2e\cdot v$, respectively; cf. proof of Theorem 8.5. The path 
$\gamma +\tilde{\gamma}$ in $\cycd{0,\bbr}{\pc^\N}$
connects $\sigma$ to $(d+2e)\cdot x_o = \deg{\sigma}\cdot x_o$. 
This shows that the degree map is an isomorphism.
\qed

We now complete the proof of Theorem 3.3, which provides a
canonical splitting of $\calz^q_\bbr$ into a product of
Eilenberg-MacLane spaces. We have a canonical homotopy equivalence
$\calz_\bbr^q \simeq \cycd{0,\bbr}{\pc^q}$, cf.  Proposition 8.1(a),
and a natural topological isomorphism
$$
\cycd{0,\bbr}{\pc^q} \cong
\pi_0(\cycd{0,\bbr}{\pc^q}) \times  \cycd{0,\bbr}{\pc^q}_o = 
\bbz \times  \cycd{0,\bbr}{\pc^q}_o 
\tag{8.18}
$$ 
given the choice of $x_o$.

\Lemma {8.8}\ {\sl
There is a canonical homotopy equivalence
$$  \cycd{0,\bbr}{\pc^q}_o \simeq 
\cycd{0}{\pc^q}_o^{av} \times \crdd{0}{\pc^q}_o.$$
}
\pf
It suffices to show that one can produce
a canonical homotopy section for the principal fibration
$\cycd{0}{\pc^q}_o^{av} \to  \cycd{0,\bbr}{\pc^q}_o \to \crdd{0}{\pc^q}_o.$
However, using (8.12) and the canonical equivalences of  diagram (8.11) 
one just needs to produce a canonical homotopy section for the
fibration
$\cycd{0}{S^{2n}}_o^{av} \to \cycd{0}{S^{2n}}_o^{fix} \to
\cycd{0}{S^n}_o\otimes \bbz_2.$

Here $S^{2n}$ is the one-point compactification of $\bbc^n$ and
$S^n$ is the sphere of  real points. 
If $i : S^n  \hookrightarrow S^{2n}$ denotes the inclusion, then let
$i_* : \cycd{0}{S^o} \to \cyfd{0}{S^{2o}}$ be the induced homomorphism.
Let $i_o : {S^n} \to \cycd{0}{S^{n}}_o$ be the usual inclusion
$i_o (x) = x - x_\infty$, where $x_\infty$ is the base point, and
let
$$
\epsilon_o : {S^n} \to \cyfd{0}{S^{2n}}_o
\tag{8.19}
$$
denote the composition $i_*\circ i_o$. 
Note that one has a commutative diagram
$$
\CD
\cycd{0}{S^n}_o @>{i_*=\epsilon_{o*}}>> \cyfd{0}{S^{2n}}_o \\
@V{q}VV @VV{p}V \\
\cycd{0}{S^n}_o\otimes \bbz_2 @>>{=}> \cyfd{0}{S^{2n}}_o /
\cyavd{0}{S^{2n}}_o , 
\endCD
\tag{8.20}
$$
where $q $ and $p$ denote the quotient maps.

Let
$f_2 : S^n \to S^n$ be the map of degree $2$, fixed in Appendix A, \S A.1.
It follows from standard properties of $H$-spaces that
$\epsilon_o \circ f_2 $ is homotopic to $2\epsilon_o$. 
On the other hand, the map $2\epsilon_o$ factors through the averaged cycles
$\cyavd{0}{S^{2n})}_o\simeq \cycd{0}{S^n\# \pr^{n-1}}_o$ and the latter space
is $(n+1)$-connected. Therefore $2\epsilon_o$ is homotopic to zero.
Now, Lemma A.1 and Corollary A.4 provide a canonical map 
$$
H : \cycd{0}{S^n}_o\otimes \bbz_2 \to \cyfd{0}{S^{2k}}_o,
$$ 
unique up to homotopy, with the property that
$$ H \circ q \simeq \epsilon_{o*}.
\tag{8.21}
$$

Let us apply Corollary A.4 once again, with 
$Y = \cycd{0}{S^n}\otimes \bbz_2$ 
and $h$ being the composition 
$S^n \xrightarrow{i_o} \cycd{0}{S^n}
\xrightarrow{q}  \cycd{0}{S^n}\otimes \bbz_2$.
It follows from (8.21) that
$ (p\circ H)\circ q \simeq  p \circ \epsilon_{o*}$,
and since $p \circ \epsilon_{o*} =q$, cf. (8.20), one concludes that
$$
(p\circ H)\circ q \simeq  q = id \circ q =  (q\circ i_o)_* .
$$
Therefore both $(p\circ H)$ and $id$ satisfy the same condition in
Corollary A.4, and since $\pi_{n+1}(\cycd{0}{S^n}\otimes \bbz_2)=0$ one
obtains $p\circ H \simeq id$. Therefore, $H$ is the desired homotopy
section of $p$.
\qed
\medskip

The combination of  Proposition 8.1(a), (8.18), Lemma 8.8 (8.16) and
Theorem 8.5 provides the proof of Theorem 3.3.
\medskip

The following facts follow from the arguments that we
have used above.

\Prop{8.9}{\sl

\roster
\item The inclusion $i : \cyav{\N}{\pc^\M} \hookrightarrow
\crl{\N}{\pc^\M}$ induces an inclusion of homotopy groups as direct 
summands;
\item The inclusion $\crl{\N}{\pc^\M} \hookrightarrow
\crl{\N+1}{\pc^{\M+1}}$ induced by the inclusion of $\pc^\M$ as a 
linear
subspace of $\pc^{\M+1}$ induces an inclusion of homotopy groups as
direct summands.
\endroster}

 \def\r{9}   \def\add{8}
\def\rc#1#2{\cz_{\bbr}^{#1}(\bbp^{#2})}
\def\rci{\cz_{\bbr}^{\infty}}
\def\ci{\cz^{\infty}}
\def\cit{\widetilde{\cz}_{\bbr}^{\infty}}
\def\ciav{\cz_{\text{av}}^{\infty}}
\def\bb#1{\overline{\b}_{#1}}
\def\la{\lambda}

\define\bzt{\bbz_2}

\subheading{\S \r.\ The ring structure}
The algebraic  join of cycles, defined for example in   [\L$_1$],  
[\L$_2$], and [\FM] is
equivariant with respect to conjugation and defines   biadditive 
pairings \newline $\rc q
n \wedge \rc {q'}{n'} \arr \rc {q+q'}{n+n'+1}$. This induces a pairing
$$
\# : \rci\wedge \rci \arr \rci ,
$$
where $\rci = \lim_{n,q\to \infty}\rc q n$, which makes $\rci$ an
$E_{\infty}$-ring space. (See \S 6.)  \ The induced map
$$
\pi_*\rci \otimes \pi_* \rci \ \arr\ \pi_* \rci
\tag{\r .1}
$$
makes $\pi_* \rci$ a graded ring.  In this section we shall compute 
this ring
and give explicit representatives for the generators.

To set the background we recall two analogous cases.  Let
$$
\ci = \lim_{n,q\to \infty}\cz^q(\bbp^n)  \qquad\text{and}\qquad
\cit =  \lim_{n,q\to \infty}\widetilde{\cz}_{\bbr}^q(\bbp^n).
$$
These are $E_{\infty}$-ring spaces as seen in [\BLLMM], and their 
homotopy groups
form graded rings.  Results from [\FM] establish an isomorphism
$$
\pi_*\ci \cong \bbz[s]
\tag{\r .2}
$$
where $s$ corresponds to the generator of $\pi_2\ci \cong \bbz$.
Results of [\Lam] show that
$$
\pi_*\cit \cong \bbz_2 [y]
\tag{\r.3}
$$
where $y$ corresponds to the generator of $\pi_1\cit \cong \bbz_2$.
The main result of this section is the following theorem  which neatly 
organizes
the additive results of \S\add.  Let $\ciav \subset \rci$ be the 
subspace defined
by taking the limits of the subgroups $\cz_{\text{av}}^q(\bbp^n) 
\subset
\cz_{\bbr}^q(\bbp^n)$ as above.   Note that the join of an averaged 
cycle with a
fixed cycle is again an averaged cycle.

\Theorem{\r.1}  {\sl  There is a ring isomorphism
$$
\pi_* \rci \ \cong\ \bbz[x,y]/(2y)
\tag{\r.4}
$$
where $x$ corresponds to the generator of $\pi_4\rci \cong \bbz$ 
and $y$
corresponds to the generator of $\pi_1\rci \cong \bbz_2$, and where 
$(2y)$
denotes the principal ideal in the polynomial ring generated by $2y$. 
Furthermore, under this isomorphism  the ideal $\pi_*\ciav \subset 
\pi_*\rci$
corresponds to the ideal
$$
\pi_*\ciav\ \cong\ (2,x)
\tag{\r.5}
$$
generated by 2 and $x$.

}

\pf
To begin we introduce a doubly indexed filtration on $\pi_*\rci$ and 
show that it is compatible with the multiplication (\r.1).  Consider the
direct  sum  $\bbc^{\infty}$ with  coordinates  $(z_0,z_1,z_2, \dots)$;
$z_j = x_j+iy_j$, and for each $n$ set 
$\bbc^{n+1} = \{z \in \bbc^{\infty} : z_j = 0 \ \text{for}\ j>n\}$.
 For each  $k$, $0\leq k\leq n$ , we consider the conjugation 
invariant
subspace 
$$
V^{n,k} \ =\ \{z \in \bbc^{n+1} \ :\    y_j=0 \  \text{for} \ j \geq k\}
\ \cong\ \bbc^k\oplus\bbr^{n+1-k}.
$$
and set 
$$
\bbp^{n,k} = \pi(V^{n,k}-\{0\})
$$
where $\pi:\bbc^{n+1}-\{0\}  \arr\bbp^n$ is the projection. 
Note that  $V^{n,k} \subset V^{n',k'}$ if $n\leq n'$ and $k\leq k'$,
and that
$$
\dim_{\bbr} \bbp^{n,k}\ =\ n+k.
$$
\noindent
{\bf Definition.} \ \ Let
$\cz^q_{\bbr}(\bbp^{n,k})\subset \rc q n$ denote the subgroup 
generated
by effective cycles $c$ for which
$$
\text{dim}_{\bbr}\left(|\widetilde c| \cap V^{n,k}\right)\ \geq \ 
\text{dim}_{\bbc}\left(|\widetilde c|\right)
\tag{\r.6}
$$
where $|\widetilde c| = \pi^{-1}(|c|) \subset \bbc^{n+1}$ denotes the
homogeneous cone corresponding  to the support  $|c|$ of $c$.
The inclusions $\cz^q_{\bbr}(\bbp^{n,k}) \subset \rc q n \subset\rci$
induce homomorphisms
$$
\pi_*\cz^q_{\bbr}(\bbp^{n,k}) \ \arr\ \pi_*\rci
\tag{\r.7}
$$

\medskip
\noindent
{\bf Observation \r.2} \ \ The image of the homomorphism (\r.7) 
remains constant 
under continuous deformations $\bbc^n$ through Real subspaces of 
$\bbc^{\infty}$

\medskip
\noindent
{\bf Observation \r.3}  \ \  The homomorphism
$$
\pi_*\cz^n_{\bbr}(\bbp^{n,k}) \arr \pi_*\rci
$$
is injective.  This follows from the results in \add.2 concerning 0-
cycles.

\medskip
\noindent
{\bf Observation \r.4}  \ \  Taking algebraic suspension by adding 
coordinates on the
left gives a commutative diagram
$$\CD
\pi_*\cz^n_{\bbr}(\bbp^{n,k}) @>{j}>>  \pi_*\rci  \\
@V{\widetilde{\csus}_*^{\ell}}VV  @V{||}V{\csus_*^{\ell}}V  \\
\pi_*\cz^n_{\bbr}(\bbp^{n+\ell,k+\ell}) @>{j_{\ell}}>>   \pi_*\rci
\endCD
\tag{\r.8}$$
where $j$ is injective by \r.3, $\csus_*^{\ell}$   is an isomorphism,
and $\widetilde{\csus}_*^{\ell}$, induced by the restriction of the 
suspension
map, is therefore also injective.

\medskip

We set 
$$
\cf^{n,k}\ \equdef\ j \pi_*\cz^n_{\bbr}(\bbp^{n,k})
$$
and note that $\cf^{n,k}$ gives a bifiltration  of $\pi_* \rci$, namely
$$
\cf^{n,k}  \subset \cf^{n',k'} \qquad 
\text{if} \ n\leq n'\ \ \text{and}\ \ k \leq k'.
$$

\medskip
\noindent
{\bf Proposition \r.5}  {\sl The homomorphism 
$\widetilde{\csus}_*^{\ell}$
in (\r.8) is an isomorphism.  Consequently, 
$$
\cf^{n,k}\ =\ \text{Im}(j_{\ell})   \qquad\text{for all $\ell > 0.$}
$$
}
\pf
In fact the map 
$\widetilde{\csus}_*^{\ell} : \cz^n_{\bbr}(\bbp^{n,k}) \to
\cz^n_{\bbr}(\bbp^{n+\ell,k+\ell})$ is a  homotopy
equivalence.  To see this one repeats the arguments of [\Lam] which 
prove that 
${\csus}_* : \cz^q_{\bbr}(\bbp^{n}) \arr \cz^q_{\bbr}(\bbp^{n+1})$ 
is a homotopy equivalence, and one notes that all  steps preserve the
subgroups $\cz^q_{\bbr}(\bbp^{n+1,*})$.  More  specifically, there are
two fundamental constructions in this proof: pulling to the normal 
cone and  ``magic fans''.  

We begin with pulling to the normal cone.  Let's introduce
homogeneous coordinates $(z,x) \in \bbc^k\oplus \bbr^{n+1-k}
 = V^{n,k} \subset \bbc^{n+1}$ and  $(\xi, z,x) \in
\bbc\oplus\bbc^k\oplus \bbr^{n+1-k}
 = V^{n+1,k+1}\subset \bbc^{n+2}$. Consider the multiplicative flow
$\varphi_t$ on $\bbp^{n+1}$ defined in homogeneous coordinates by
$\varphi_t(\xi, z,x) = (t\xi, z,x)$ for $t \in \bbr^+$. This flow
induces ``pulling to the  normal cone'' in [\Lam]. It evidently
preserves condition (\r6) above, and therefore preserves the
subgroup of algebraic cycles $\cz^q_{\bbr}(\bbp^{n+1,k+1})$.

In the ``magic fan'' construction one adds a new coordinate giving
$(\eta, \xi, z,x) \in \bbc\oplus\bbc\oplus\bbc^k\oplus \bbr^{n+1-k}
= \bbc \oplus  V^{n+1,k+1}$.
To each homogeneous polynomial $f(\eta, \xi, z,x)$ with real 
coefficients one
constructs a transformation  $\Phi_f : \cz^q_{\bbr}(\bbp^{n+1}) \to 
\cz^q_{\bbr}(\bbp^{n+1})$ by setting $\Phi_f(c) = 
(\pi_1)_*\left(\pi_0^*c\bullet D_f\right)$ where $D_f$ is the divisor
of $f$, and $\pi_0, \pi_1$ are  projections
$\bbp_{\bbc}^{n+2}\ \cdot\cdot\cdot\cdot  >\ \bbp_{\bbc}^{n+1}$ 
with vertices $(1,0,0,\dots,0)$ and $(1,1,0,\dots,0)$ respectively. 
We need to  check that this construction preserves the subgroups 
$\cz^q_{\bbr}(\bbp^{n+1,k+1})$.  For this let $\widetilde Y
\subset \bbc^{n+2}$ denote the homogeneous cone of a projective
variety $Y\subset \bbp^{n+1}$.  It will suffice to show that
$$
\dim_{\bbr}\left(\widetilde Y \cap V^{n+1,k+1}\right)
\ \leq\ 
\dim_{\bbr}\left(\widetilde{\left(\Phi_f Y\right)} \cap
V^{n+1,k+1}\right) 
\tag{\r.9}
$$
To see this consider a point $a \in Y$ with homogeneous coordinates
$(\xi, z,x) \in \widetilde Y \cap V^{n+1,k+1}$.    Let $\eta_1, \dots
, \eta_d$ be the zeros of  the polynomial $q(t) = f(t, \xi, z,x)$.
Then $\pi_0^{-1}(a) \cap D_f $ consists of the $d$ points with
homogeneous coordinates $(\eta_j, \xi, z,x)$, \ $j=1,\dots, d$, and
$\pi_1\left(\pi_0^{-1}(a) \cap D_f \right)$ is the union of the $d$
points with homogeneous coordinates 
$(\eta_j-\xi, z,x)$, \ $j=1,\dots, d$.  Note that each of these
points again lies in $V^{n+1,k+1}$.  It follows that condition (9.9)
holds as claimed.   Therefore the arguments of [\Lam]   apply
without change to  show that 
$\widetilde{\csus}$ is a homotopy equivalence, and we are done.  
\qed
\medskip

The  analogues of \r.2 -- \r.5   apply also to the averaged cycles:

$$
\cz^q_{\text{av}}(\bbp^{n,k}) =  \cz^q_{\bbr}(\bbp^{n,k}) \cap
\cz^q_{\text{av}}(\bbp^{n}) 
$$
One obtains a bifiltration $\cf^{n,k}_{\text{av}}$ of $\pi_*\ciav$
where the $n$-filtration agrees with that of Theorem \add.5.  Under 
the 
isomorphism $\pi_*\ciav \cong \widetilde{H}_*(X^\infty ; \bbz)$  
deduced in \S
\add \ (See \add.3), consider the classes
$$
\theta_{n,k} \ =\ [\bbp^{n,k}/\bbz_2] \in 
\widetilde{H}_{n+k}(X^\infty ; \bbz)
\cong \pi_{n+k}\ciav
\qquad\text{for }\ \ 0<k\leq n.
\tag{\r.10}
$$
Note that the  intersection of $\bbp^{n,k}$ with the affine
coordinate chart $\bbp_{\bbc}^{n}-\bbp_{\bbc}^{n-1}
 = \bbc^n = \bbr^n \oplus i\bbr^n$ is exactly $\bbr^n \oplus i\bbr^k$.
Therefore, the image of $\theta_{n,k}$ in the homology of $X^n/X^{n-
1} =
S^{2n}/\bbz_2$ is the class $S^{n}\# \bbp_{\bbr}^{k-1}$.  Hence this 
image
generates the group $\Iav{n}{k}$ in the bigrading established in 
Theorem \add.5.

It follows that
   
$$
\cf^{n,k}_{\text{av}} \ =\ \text{span}_{\bbz}\{\theta_{n',k'}\, :\, n' 
\leq n
\ \ \text{and}\ \ k'\leq k\},
 $$
and furthermore that under the injection $\pi_*\ciav 
\hookrightarrow
\pi_*\rci$ we have that 
$$
\cf^{n,k}_{\text{av}}\ =\ \cf^{n,k} \cap \pi_*\ciav.
\tag{\r.11}
$$
In particular, we deduce that for $0< k \leq n$
$$\aligned
\GRav{n}{k}  &\equdef \cf^{n,k}_{\text{av}}\cap 
\pi_{n+k} \cz^{\infty}_{\text{av}}
= \langle \theta_{n,k}\rangle  \\
&\cong  \cf^{n,k}_{\text{av}}/
\cf^{n-1,k}_{\text{av}} \oplus \cf^{n,k-1}_{\text{av}} 
\endaligned$$
and $\GRav{n}{0} = \{0\}$.  For $0\leq k\leq n$ set
$\GR{n}{k}  \equiv \cf^{n,k}\cap \pi_{n+k} \cz^{\infty}_{\bbr}$.
From Theorem \add.2 we know that $\GR{n}{0} \cong 
\pi_n\{\cz_0(\bbp-{\bbc}^n/\bbp-{\bbc}^{n-1})_o\}^{fix} = \bbz_2$.
Let $\theta_{n,0}$ denote the generator of $\GR{n}{0}$. (an explicit
representative will be given later).  Then we have that for $0 \leq k 
\leq n$
$$\aligned
\GR{n}{k}  &\equdef \cf^{n,k}\cap 
\pi_{n+k} \cz^{\infty}_{\bbr}
= \langle \theta_{n,k}\rangle  \\
&\cong  \cf^{n,k}/
\cf^{n-1,k} \oplus \cf^{n,k-1} 
\endaligned$$
It is useful to picture the graded peices $\GR{n}{k}$
on the $(n,k)$-coordinate grid.

\medskip

 $$
\matrix     k\uparrow
     & \   & \   & \   & \   & \   & \   & \   & \   & \   & \   &\dots \\  
 \   & \   & \   & \   & \   & \   & \   & \   & \   & \   &\bbz &\dots \\
 \   & \   & \   & \   & \   & \   & \   & \   & \   & 0   & 0   &\dots \\
 \   & \   & \   & \   & \   & \   & \   & \   &\bbz &\bzt &\bzt &\dots 
\\
 \   & \   & \   & \   & \   & \   & \   & 0   & 0   & 0   & 0   &\dots \\
 \   & \   & \   & \   & \   & \   &\bbz &\bzt &\bzt &\bzt &\bzt 
&\dots \\
 \   & \   & \   & \   & \   & 0   & 0   & 0   & 0   & 0   & 0   &\dots \\
 \   & \   & \   & \   &\bbz &\bzt &\bzt &\bzt &\bzt &\bzt &\bzt 
&\dots \\
 \   & \   & \   & 0   & 0   & 0   & 0   & 0   & 0   & 0   & 0   &\dots \\
 \   & \   &\bbz &\bzt &\bzt &\bzt &\bzt &\bzt &\bzt &\bzt &\bzt 
&\dots \\
 \   & 0   & 0   & 0   & 0   & 0   & 0   & 0   & 0   & 0   & 0   &\dots \\ 
 \bbz&\bzt &\bzt &\bzt &\bzt &\bzt &\bzt &\bzt &\bzt &\bzt &\bzt 
&\dots& \underset{n}\to \longrightarrow
\endmatrix
$$ 

 \medskip

The graded subgroup $\GRav{n}{k}\subset \GR{n}{k}$
looks like this:

\medskip

$$
\matrix     k\uparrow  
     & \   & \   & \   & \   & \   & \   & \   & \   & \   & \   &\dots \\
 \   & \   & \   & \   & \   & \   & \   & \   & \   & \   &\bbz &\dots \\
 \   & \   & \   & \   & \   & \   & \   & \   & \   & 0   & 0   &\dots \\
 \   & \   & \   & \   & \   & \   & \   & \   &\bbz &\bzt &\bzt &\dots 
\\
 \   & \   & \   & \   & \   & \   & \   & 0   & 0   & 0   & 0   &\dots \\
 \   & \   & \   & \   & \   & \   &\bbz &\bzt &\bzt &\bzt &\bzt 
&\dots \\
 \   & \   & \   & \   & \   & 0   & 0   & 0   & 0   & 0   & 0   &\dots \\
 \   & \   & \   & \   &\bbz &\bzt &\bzt &\bzt &\bzt &\bzt &\bzt 
&\dots \\
 \   & \   & \   & 0   & 0   & 0   & 0   & 0   & 0   & 0   & 0   &\dots \\
 \   & \   &\bbz &\bzt &\bzt &\bzt &\bzt &\bzt &\bzt &\bzt &\bzt 
&\dots \\
 \   & 0   & 0   & 0   & 0   & 0   & 0   & 0   & 0   & 0   & 0   &\dots \\ 
 2\bbz& 0  & 0   & 0   & 0   & 0   & 0   & 0   & 0   & 0   & 0   &\dots&
\underset{n}\to \longrightarrow
\endmatrix
$$

\medskip

Note that $\GR{n}{k} = \GRav{n}{k}$  for  $k > 0$ and 
$\GRtil{*}{*} \equiv \GR{*}{*}/\GRav{*}{*}$
is simply:

\medskip

$$
\matrix     k\uparrow  
 \   & \   & \   & \   & \   & 0   & 0   & 0   & 0   & 0   & 0   &\dots \\
 \   & \   & \   & \   & 0   & 0   & 0   & 0   & 0   & 0   & 0   &\dots \\
 \   & \   & \   & 0   & 0   & 0   & 0   & 0   & 0   & 0   & 0   &\dots \\
 \   & \   & 0   & 0   & 0   & 0   & 0   & 0   & 0   & 0   & 0   &\dots \\
 \   & 0   & 0   & 0   & 0   & 0   & 0   & 0   & 0   & 0   & 0   &\dots \\ 
 \bzt&\bzt &\bzt &\bzt &\bzt &\bzt &\bzt &\bzt &\bzt &\bzt &\bzt 
&\dots&
\underset{n}\to \longrightarrow
\endmatrix
$$

\medskip

Our main observation here is the following.

\Prop {\r.6} {\sl The filtrations $\cf^{n,k}$ and
$\cf^{n,k}_{\text{av}}$ are compatible with the join pairing, i.e.,}
$$
\#\left(\cf^{n,k}\otimes\cf^{n',k'}\right) \subset \cf^{n+n',k+k'}
\qquad\text{and}\qquad
\#\left(\cf^{n,k}_{\text{av}}\otimes\cf^{n',k'}_{\text{av}}\right) 
\subset \cf^{n+n',k+k'}_{\text{av}}
$$
 
\pf
By Observation \r.2 the join gives a well defined homomorphism
$$
\pi_*\cz^{n}_{\bbr}(\bbp^{n,k}) \otimes 
\pi_*\cz^{n'}_{\bbr}(\bbp^{n',k'}) 
\ \arr\ \pi_*\cz^{n+n'}_{\bbr}(\bbp^{n+n'+1,k+k'+1}).
$$
Pushing into $\pi_*\rci$ and applying \r.5 gives the result for 
$\cf^{*,*}$.
The argument for $\cf^{*,*}_{\text{av}}$ is similar. \qed

\Lemma{\r.7}  {\sl  The generators 
$\theta_{2,2} \in \pi_4\rci \cong \pi_4 \ciav \cong \bbz$ and 
$\theta_{1,0} \in \pi_1\rci$ have the property that
$$
\theta_{2,2} \in \cf^{2,2} \qquad\text{and}\qquad \theta_{1,0} \in 
\cf^{1,0}
$$
}

\pf
Recall $X^q \equiv \bbp_{\bbc}^q/\bbz_2$. (See Cor. 8.4.)\ 
Under the isomorphisms $\pi_4 \cz^q_{\text{av}} \cong \pi_4 
\cz_0(X^q) \cong
H_4(X^q;\bbz) \cong \bbz$, the generator $\theta_{2,2}$  corresponds 
to the class
$[X^2] \in H_4(X^q;\bbz)$ for any $q\geq 2$.  Thus $\theta_{2,2} \in 
\cf^{2,2}$.
For the second assertion recall that the generator of $\pi_1\rci$ is 
given by
the map $S^1 \arr \cz_0(\bbp^1_{\bbr}) \subset 
\cz_0(\bbp^1_{\bbc})$ sending $t
\mapsto t-t_0$ for $t \in \bbp^1_{\bbr} \cong S^1$ (where $t_0 \in 
\bbp^1_{\bbr}$
is a base point).  \qed

\Prop {\r.8} {\sl

\hskip 1.2in A)\ \ $\theta_{2,2}^k = \theta_{2k,2k}\ \ \text{in}\ \
\pi_{4k}\ciav.$

\hskip 1.2in B)\ \ $\theta_{2,2}^k\cdot \theta_{1,0}^\l \neq 0 \ \ 
\text{in}\ \
\pi_{4k+\l}\ciav.$ }
\bigskip
\noindent
{\bf Proposition \r.8 $\Rightarrow$ Theorem \r.1.} \ \ From results 
in \S\add\  we
have an exact sequence
$$
0\arr \pi_* \ciav \arr \pi_* \rci @>{p_*}>> \pi_*\cit \arr 0
\tag{\r.12}
$$
and we know that the elements $p_*(\theta_{1,0}^{\ell}), \ \ \ell 
\geq 1$, give
an additive basis of $\pi_*\cit$.  In fact, $p_*$ is a ring 
homomorphism and 
$p_*(\theta_{1,0})$ corresponds to $y$ in (\r.3). 

We have seen that under the identification
$$
\pi_*\ciav \ \cong \ H_*(X^{\infty};\bbz)
$$
the bifiltration of $\pi_*\ciav$ corresponds to the bifiltration of 
$H_*(X^{\infty};\bbz)$ induced by the family of subspaces 
$\bp^{n,k}/\bbz_2$ in $X^{\infty}$.  From the results of \S\add \ we 
have   
 that
$$
\cf_{\text{av}}^{ *,\text{odd}} = 0
$$
and 
$$
\cf_{\text{av}}^{2k+\ell,2k}\cap \pi_{4k+\ell}\ciav \ =\ 
\cases &\bbz\cdot \theta_{2k,2k}\ \ \text{if}\ \  \ell = 0  \\
&\bbz_2\cdot \theta_{2k+\ell,2k}\ \ \text{if} \ \ \ell > 0.\endcases 
$$
Since 
$$
\theta_{2,2}^k\cdot \theta_{1,0}^\ell \in 
\cf_{\text{av}}^{2k+\ell,2k}\cap \pi_{4k+\ell}\ciav,
$$
the conclusion of Proposition \r.8 implies that 
$\theta_{2,2}^k\cdot \theta_{1,0}^\ell = \theta_{2k+\ell,2k}$.
\qed

\noindent
{\bf Proof of Proposition \r.8} \ \ We begin by constructing explicit
representatives of $\theta_{1,0}$ and $\theta_{2,2}$.

Consider $S^1 = \bbp^1_{\bbr} \subset \bbp^1_{\bbc}$, the fixed-
point set,
and choose a base point $t_0 \in S^1$.
Define
$$
\a : S^1 \arr \rc 1 1 
\tag{\r.13}
$$
by 
$$
\a(t) = t-t_0.
$$
This map clearly has filtration level (1,0) since the image is 
supported in 
$\bbp^1_{\bbr}$.  One sees directly that under the projection
$p_* :\rc 1 1 \arr \widetilde{\cz}^1_{\bbr}(\bbp^1_{\bbc}) = 
\cz_0(\bbp^1_{\bbr})\otimes \bbz_2$ the map $p_*\circ \a$ 
represents the
generator of $\pi_1(\cz_0(\bbp^1_{\bbr})) \cong H_1(\bbp^1_{\bbr}; 
\bbz_2) \cong
\bbz_2$.  Hence $\a$ represents  the non-zero class $\theta_{1,0}$ in 
$\pi_1\rci
= \bbz_2$.

Consider $\bbp^1_{\bbr}\subset \bbp^1_{\bbc}$ as the ``equator'' 
and let $D^2
\subset \bbp^1_{\bbc}$ be the ``upper hemisphere'', (so 
$\bbp^1_{\bbc} = D^2
\cup\overline{D^2}$  where $\overline{(\cdot)}$ is the conjugation 
map).  For each
$n\geq 1$ we define a map
$$
\b_n : (D^2)^n \ \arr\ \cz^n(\bbp_{\bbc}^{2n-1})
\tag{\r.14}
$$
by
$$
\b_n(t_1,\dots,t_n) \ = \ (t_1 - \bar{t}_1)\#\dots \#(t_n - \bar{t}_n).
$$
Note that $\b_n(t_1,\dots,t_n) = 0$ if $t_j \in \partial D^2 =  
\bbp_{\bbr}^1$
for {\bf any} $j$.  Thus $\b_n$ descends to a map
$$
\bb n : S^2\wedge \dots \wedge S^2 = S^{2n} \arr 
\cz^n(\bbp_{\bbc}^{2n-1}).
$$
Note that $\overline{\b_1(t)} = -\b_1(t)$, that is $\b_1$ maps into 
{\sl anti-averaged} cycles.  The join of two anti-averaged cycles in an 
averaged
cycle.  Since $\b_n(t_1,\dots,t_n) = \b_1(t_1) \#\dots\#\b_n(t_n)$ 
we see that
$$
\bb n : S^{2n} \arr \cz^n_{\text{av}}(\bbp_{\bbc}^{2n-1})
$$
{\sl whenever $n$ is even.}

From Observation \r.5 we see that the class  $[\bb 2]$ of $\bb 2 : S^4 
\to
\cz^2_{\text{av}}(\bbp^3) = \cz^2_{\text{av}}(\bbp^{3,3})$ in $\pi_4 
\ciav$ has
filtration level (2,2).
It follows that the class $[\bb{2n}] = [\bb 2]^n \in \pi_{4n}\ciav$ has
filtration level $(2n,2n)$.

\Lemma{\r.9}  {\sl The map $\bb n : S^{2n} \arr 
\cz^n(\bbp_{\bbc}^{2n-1})$
represents the generator of $\pi_{2n} \cz^{\infty} \cong \bbz$.}

\pf
Recall that $\b_1 : (D^2,\partial D^2) \arr (\cz_0(\bbp_{\bbc}^{1}, 0)$ 
is given
by $\b_1(t) = t - \bar t$.  To compute the class of $\bb 1$ in 
$\pi_2\cz_0(\bbp_{\bbc}^{1}) = H_2(\bbp_{\bbc}^{1}; \bbz)$ we take 
the graph
$\Gamma(\bb 1)$ in $S^2\times \bbp_{\bbc}^{1}$ and push it 
forward to
$\bbp_{\bbc}^{1}$  (cf. [\FL]).  Now $\Gamma(\bb 1)$ is an oriented 
cycle which
is the union of two oriented disks $\Gamma_0\cup \Gamma_1$.  The 
disk
$\Gamma_0$ is the graph of the identity map $D^2 \to D^2_+$ on the 
upper
hemisphere with the canonical orientation from $D^2$; the disk
$\Gamma_1$ is the graph of the conjugation map $D^2 \to D^2_-$ 
from the upper
to the lower hemisphere with the   orientation opposite the one given 
by
 $D^2$ due to the minus sign in $\bb 1$.  Note that 
$\Gamma_0\cup \Gamma_1 \subset S^2\times \bbp_{\bbc}^{1}$ is 
an oriented
2-sphere homeomorphic to $\bbp_{\bbc}^{1}$ under projection to the 
second factor.

\vskip 1.5in

This shows that $\bb 1$ represents the generator $s$ of 
$\pi_2\cz^{\infty}$.
It follows that $\bb n = \bb 1 \#\dots \# \bb 1$ represents 
$s^n \in \pi_{2n}\cz^{\infty}$, which is the generator by \r.2. \qed

\Lemma {\r.10} {\sl For any $N \geq 2n$, the homomorphism
$I_* : \pi_{4n}\cz^{2n}_{\text{av}}(\bbp_{\bbc}^N) \arr \pi_{4n} 
\cz^{\infty}$,
induced by the inclusion $i:\cz^{2n}_{\text{av}}(\bbp_{\bbc}^N)
\hookrightarrow \cz^{\infty}$, is an isomorphism.}

\pf
By the Algebraic Suspension Theorem [\L$_1$], [\LLM$_2$] it 
suffices to consider the map
of 0-cycles $\cz^{2n}_{\text{av}}(\bbp_{\bbc}^{2n}) \arr 
\cz^{2n}(\bbp_{\bbc}^{2n})$.  Note that the composition
$$
\cz^{2n}_{\text{av}}(\bbp_{\bbc}^{2n}) @>{i}>>
\cz^{2n}(\bbp_{\bbc}^{2n})     @>{\text{av}}>>
\cz^{2n}_{\text{av}}(\bbp_{\bbc}^{2n})
$$
where $\text{av}(x) = x + \bar x$, is multiplication by 2, and so 
therefore is
the composition
$$\CD
\pi_{4n}\cz^{2n}_{\text{av}}(\bbp_{\bbc}^{2n}) @>{i_*}>>
\pi_{4n}\cz^{2n}(\bbp_{\bbc}^{2n})     @>{\text{av}_*}>>
\pi_{4n}\cz^{2n}_{\text{av}}(\bbp_{\bbc}^{2n}) \\
||  @.  || @. || \\
\bbz @. \bbz @.\bbz 
\endCD$$
On the other hand the homomorphism av$_*$ can be identified with 
the
homomorphism $\rho_* : H_{4n}(\bbp_{\bbc}^{2n}; \bbz) \arr
H_{4n}(\bbp_{\bbc}^{2n}/\bbz_2; \bbz) $  where $\rho : 
\bbp_{\bbc}^{2n}
\arr \bbp_{\bbc}^{2n}/\bbz_2$ is the quotient map.  This map 
clearly sends the
fundamental class $[\bbp_{\bbc}^{2n}]$ to 
$2[\bbp_{\bbc}^{2n}/\bbz_2]$.
Hence av$_* = 2$ and so $i_*$ must be an isomorphism. \qed

\Cor {\r.11} {\sl The map $\bb 2$ represents the generator 
$\theta_{2,2}$ of 
$\pi_4\cz_{\text{av}}^2(\bbp_{\bbc}^3)$.  Furthermore, for all 
$n\geq 1$ one has
that $\theta_{2,2}^n = \theta_{n,n}$.
}

\pf  
The first assertion follows immediately from \r.9 and \r.10.  For the 
second we
recall that $[\bb{2n}] = [\bb 2]^n$ in $\pi_* \ciav$, and use (\r.2).  
\qed

\medskip

This establishes part A) of Proposition \r.8.  For part B) we invoke 
the
following.  Consider a continuous map $f:S^m \to
\cz_{\text{av}}^q(\bbp_{\bbc}^{N})$ and let $\rho : \bbp_{\bbc}^{N}
\to X^{N}$ be the projection.  Then $\rho_* f(x) = 2\bar f (x)$ where
$\bar f:S^m \to\cz^q(X^{N})$ is a continuous map into cycles on 
$X^{N}$.
Assume that $f$ is well enough behaved to have a graph 
$\Gamma_f$ in 
$S^m\times\bbp_{\bbc}^{N}$ (cf. [\FL]).  Then $(1\times 
\rho_*)[\Gamma_f]
=2[\Gamma_{\bar f}]$ where $\Gamma_{\bar f}$ is a cycle on 
$S^m\times X^{N}$ which we will call the graph of $\bar f$.  Let 
$\text{pr} : S^m\times X^{N} \to  X^{N}$ be projection.

\Lemma {\r.12} {\sl  If $[{\roman pr}_*\Gamma_{\bar f} ] \neq 0$
in $H_*(X^{N}; \bbz)$, then $[f] \neq 0$ in 
$\pi_m\cz_{\text{av}}^q(\bbp_{\bbc}^{N})$.
}

\pf
Suppose $f:\partial D^{m+1}\to \cz_{\text{av}}^q(\bbp_{\bbc}^{N})$
extends to a continuous map
$F: D^{m+1}\to \cz_{\text{av}}^q(\bbp_{\bbc}^{N})$ which we may 
assume 
to have a graph.  Then the integral chain pr$_* \Gamma_{\bar F}$ 
has boundary
pr$_* \Gamma_{\bar f}$ in $X^{N}$.  \qed

\medskip

To detect homology classes in $X^{N}$ we will use the following.

\Lemma {\r.13} {\sl  Let $Z\subset X^{N}$ be an integral cycle of 
codimension
$\ell < N$ defined by the oriented regular set of a real analytic 
subvariety.
  Let $\bbp_{\bbr}^{N}\subset X^{N}=\bbp_{\bbc}^{N}/\bbz_2$
denote the singular set of $X^{N}$.  Suppose there exists a compact 
oriented
submanifold $Y^{\ell} \hookrightarrow X^{N}-\bbp_{\bbr}^{N}$ of 
dimension
$\ell$ which meets the regular set of $Z$ transversely in one point.  
Then $[Z]
\neq 0$ in $H_{2N-\ell}(X^{N};\bbz)$.
 }

\pf  Let $M = X^{N}-(\bbp_{\bbr}^{N})_\epsilon$ where 
$(\bbp_{\bbr}^{N})_\epsilon$ is a tubular neighborhood of 
$\bbp_{\bbr}^{N}$ whose closure does not meet $Y^{\ell}$.  
Note that $M$ is a smooth compact oriented manifold with boundary.
The restriction  $Z_\epsilon \equiv Z\cap M$ defines an integral cycle
of codimension $\ell$ on $(M, \partial M)$, and $Y^\ell$ defines a 
cycle of
dimension $\ell$ on $M$.  The intersection hypothesis implies that
$[Z_\epsilon] \neq 0$ in $H_{2N-\ell}(M, \partial M) \cong
H_{2N-\ell}(X^{N},\bbp_{\bbr}^{N})$.  In the long exact sequence for
the pair $(X^{N},\bbp_{\bbr}^{N})$ we have  
$$
H_j(\bbp_{\bbr}^{N}) \arr H_j(X^{N}) @>{r}>> 
H_j(X^{N},\bbp_{\bbr}^{N}),
$$
and it is clear from the construction that  $r([Z]) = [Z_\epsilon]$.  \qed

\medskip

To complete the proof of Proposition \r.8 B) we consider the map
$$
f \equiv {\bb 2 }^n  \a^\ell : S^{4n+2\ell} \arr
\cz^{2n+\ell}(\bbp_{\bbc}^{4n+\ell-1}) 
$$
where we may assume that $\ell$ is odd.  This map can be 
coordinatized as
follows.  Choose affine coordinates $x_1,\dots,x_{2n}, 
y_1,\dots,y_\ell$ on 
$\bbp_{\bbc}^{1}\times \dots \times \bbp_{\bbc}^{1}$ ($(2n+\ell)$-
times)
and restrict them to
$$
\text{Im} (x_i) \geq 0 \qquad\text{and} \qquad
\text{Im} (y_j) = 0 \qquad\text{for all}\ \ i,j.
$$
Let $\la_{x_i} = \text{span}(1,x_i)$ and 
$\la_{y_j} = \text{span}(1,y_j)$ in $\bbc^2$.  Then 
$$
f(x,y) = (\la_{x_1}-\la_{\overline{x_{1}}}) 
\#\dots \#
(\la_{x_{2n}}-\la_{\overline{x_{2n}}}) \#
(\la_{y_1}-\la_{0}) 
\#\dots \#
(\la_{y_{\l}}-\la_{0}).
$$
Note that $f = \text{av} \circ \bar f$ where 
$$
\bar f(x,y) = \la_{x_1}\# (\la_{x_2}-\la_{\overline{x_{2}}}) 
\#\dots \#
(\la_{x_{2n}}-\la_{\overline{x_{2n}}}) \#
(\la_{y_1}-\la_{0}) 
\#\dots \#
(\la_{y_{\ell}}-\la_{0}).
$$
Consider the affine chart $\bbc^{4n + 2\ell - 1}$ 
on $\bbp(\bbc^2\oplus\dots\oplus\bbc^2) = \bbp(\bbc^{4n + 2\ell })$ 
given
be setting the first coordinate equal to 1.  Let $\bbc^{N}/\bbz_2$ be 
its image
in $X^{N}$ where $N = 4n + 2\ell -1$.  Then $Z = 
\text{pr}_*(\Gamma_{\bar f})$
is an analytic cycle whose image in $\bbc^{N}/\bbz_2$ is described 
as follows. 
Note that $\bar f$ expands into $2^{2n-1+\ell}$ factors.  However 
those factors
containing the constant $\la_0$ project to 0 for dimension reasons.  
Thus it
suffices to consider the remaining $2^{2n-1}$ factors.  Upstairs in 
$\bbc^{N}$ they are subsets of the form
$$\aligned
Z_{\pm  \pm \dots\pm}\ & =\
  \bigcup \biggl\{x_1 + 
\text{span} \biggl(
(0,1,x_2,0,...,0), (0,0,0,1,x_3,0,...,0), \dots , (0,\dots,0,1,x_{2n},0,...,0),
   \\
& \qquad\qquad\qquad\qquad\qquad\qquad\qquad
(0,\dots,0,1,y_1,0,...,0),\dots,(0,\dots,0,1,y_\ell)\biggr)\biggr\}
\endaligned$$
where the union is over all $x,y$ with 
$\text{Im}(y_i) =0$ for all $i$ and $\pm\text{Im}(x_j) \geq 0$ 
depending on the choice of $+$ or $-$ in the $j^{\text{th}}$ subscript 
of 
$Z_{\pm \pm \dots\pm}$.
These sets can be rewritten as 
$$\aligned
Z_{\pm \pm \dots\pm}\ &=\ \bigl\{(z_2,...,z_{4n},w_1,...,w_{2\ell}) 
\,:\,
\text{Im}(z_2) \geq 0,\\
&\qquad \pm \text{Im}(\overline{z_{2j-1}}z_{2j}) \geq 0 \ \ \forall 
j>2, 
\ \ \text{and}\ \ \text{Im}(\overline{w_{2i-1}}w_{2i}) = 0\ \  \forall 
i\bigr\}
\endaligned $$
 The union of these, with orientations adjusted for signs, is the 
oriented
semi-analytic set
$$
\widetilde Z \ =\ \bigl\{(z,w) \in \bbc^{4n-1}\times \bbc^{\ell}\,:\,
\text{Im}(z_2) \geq 0  
\ \ \text{and}\ \ \text{Im}(\overline{w_{2i-1}}w_{2i}) = 0\ \  \forall 
i  \bigr\}
$$
Thus, in this coordinate chart $\bbc^{N}/\bbz_2$ our total cycle
$\Gamma_{\bar f}  \subset  X^{N}$
is exactly the {\bf reduced} image of the real analytic variety defined 
by 
the equations $\text{Im}(\overline{w_{2i-1}}w_{2i}) = 0$, i.e., 
$$
\Gamma_{\bar f}\cap (\bbc^{N}/\bbz_2)\ =\ 
{\ssize\frac 1 2}\rho_*\bigl\{(z,w) \,:\, 
\text{Im}(\overline{w_{2i-1}}w_{2i}) = 0\ \  \forall i\bigr\} 
$$
where $\rho:\bbc^{N} \to\bbc^{N}/\bbz_2$ is the projection.

Consider now the sphere
$$
\widetilde Y \ =\ \biggl\{ (0,0,\dots,0,1,it_0, 1,it_1, \dots, 1,it_\ell)
\,:\, t_i \in \bbr\ \ \forall i\ \ \text{and}\ \ 
\sum_i t_i^2 = 1 \biggl\}
$$
and let $Y = \rho(\widetilde Y) \cong \bbp_{\bbr}^\ell$ be its 
reduced image in
$\bbc^{N}/\bbz_2\subset X^{N}$.  Note that $Y$ misses the singular 
set
$\rho(\bbr^{N})$ and $Y$ meets $\Gamma_{\bar f}$ in exactly one 
point, namely
the conjugate pair corresponding to $t_1= \dots = t_\ell = 0$ and $t_0 
= \pm 1$.
One easily checks that this is a regular point of $\widetilde Y$.
This completes the proof.   \qed\ $\square$

\define\ov{\overline}
\subheading{\S A.\ \ Appendix: Splittings and Eilenberg-MacLane spaces}
\bigskip

\noindent{\bf A.1.\ \ Models for Eilenberg-MacLane spaces}
\medskip

In our discussion, the preferred model for the Eilenberg-MacLane space
$K(\bbz,n)$ is $\cycd{0}{S^n}_o$, the connected component of $0$ in the
topological abelian group $\cycd{0}{S^n}$. If $N$ is a finitely
generated abelian group, then $\cycd{0}{S^n}_o\otimes_\bbz N$ is our
model for $K(N,n)$. Note that this model for $K(N, n)$ is a
topological $R$-module, and in particular, this model for
$K(\bbz_p,n)$ is a  {\it $p$-torsion group}.

Now, for each $p$ fix a map $f_p : S^n \to S^n$ of degree $p$, and define 
$$
M(\bbz_p,n) \equdef D^{n+1}\cup_{f_p} S^n.
$$
This is a Moore space satisfying
$$
\widetilde{H}_j(M(\bbz_p,n);\bbz) \cong
    \cases
        \bbz_p & \text{ if } j = k \\
        0 & \text{ otherwise. }
    \endcases
$$
It follows that we have yet another model for $K(\bbz_p,n)$, namely,
the {\it torsion free} abelian topological group
$\cycd{0}{M(\bbz_p,n)}_o$. We need to establish a few properties of
$M(\bbz_p,n)$ and of $\cycd{0}{M(\bbz_p,n)}_o.$

Consider the canonical inclusion $\imath : S^n \hookrightarrow
M(\bbz_p,n)$, and let $F_p : D^{n+1} \to M(\bbz_p,n)$ be the canonical
map which induces the relative homeomorphism 
$F_p : (D^{n+1},\; S^n) \to (M(\bbz_p,n),\; S^n)$ and satisfies
${F_p}_{|\partial D^{n+1}} = f_p$. 

The following result is rather standard.

\Lemma {A.1}\ {\sl Given a map $h: S^n \to Y$ such that $h\circ f_p$ is
homotopic to zero, then there is an extension of $h$ to $M(\bbz_p,n)$.
In other words, there is an $\overline{h} : M(\bbz_p,n) \to Y$ such that
$\overline{h} \circ \imath = h$. Furthermore, if $\pi_{n+1}(Y) = 0$,
then the extension is unique up to homotopy.
}

\Cor{A.2}\ {\sl Let $h : S^n \to Y$ be as in the Proposition, and
assume that $Y$ is an abelian topological group.
\roster
\item If $h$ sends the base-point $x_\infty\in S^n$ to $0\in Y$, then one
has a commutative diagram
$$
\CD
S^n @>{\imath}>> M(\bbz_p,n) @>{\overline{h}}>> Y \\
@V{j_S}VV @VV{j_M}V @VV{=}V \\
\cycd{0}{S^n}_o @>>{\imath_*}> \cycd{0}{M(\bbz_p,n)} @>>{\overline{h}_*}> Y,
\endCD
$$
where $j_S$ and $j_M$ are natural inclusions, and $\imath_*$ and
$\overline{h}_*$ are the group homomorphisms induced by $\imath$ and
$\overline{h}$, respectively. 
\item If $\pi_{n+1}(Y)=0$ then any homomorphism $\phi :
\cycd{0}{M(\bbz_p,n)}_o \to Y$, with the property 
$\phi\circ \imath_*\circ j_S = h$, is
homotopic to $\overline{h}_*$ through group homomorphisms.
\endroster
}
\pf The first assertion is a direct consequence of the fact that 
$Y$ is an abelian topological group and from universal properties of the free
abelian group on $M(\bbz_p,n)$.
To prove the second assertion, consider the map 
$ \phi\circ j_M : M(\bbz_p, n) \to Y$. 
Since 
$(\phi\circ j_M) \circ \imath = \phi  \circ \imath_* \circ j_S = h$, 
one concludes from the proposition that 
$\phi\circ j_M $ is homotopic to $\overline{h}$, and hence
$\phi = (\phi\circ j_M)_*$ is homotopic to $\overline{h}_*$ through
homomorphisms. 
\qed

\medskip
\Cor{A.3}\ {\sl Let $q : \cycd{0}{S^n}_o \to  \cycd{0}{S^n}_o\otimes
\bbz_p$ denote the quotient map.  Then there is a canonical 
homotopy equivalence
$\Psi \ : \ \cycd{0}{M(\bbz_p,n)}_o \to \cycd{0}{S^n}_o\otimes \bbz_p$ 
satisfying
\roster
\item $\Psi$ is a group homomorphism;
\item $\Psi\circ \imath_* = q.$
\endroster
Furthermore, any $\phi$ satisfying the above properties
is homotopic to $\Psi$ through homomorphisms.
}
\pf
Just observe that the composition $(q\circ j_S)\circ f$ is homotopic to
multiplication by $p$ in the homotopy group 
$\pi_n(\cycd{0}{S^n}_o\otimes \bbz_p)$, 
and hence $q\circ j_S$ satisfies the hypothesis of the proposition. 
It is easy to see that the  homomorphism resulting from 
Corollary A.2 induces an isomorphism
of $n$-th homotopy groups, and since $q = (q\circ j_S)_*$, the corollary
now follows.
\qed
\Cor{A.4}\ {\sl Given an abelian topological group $Y$ and a map $h: S^n
\to Y$ such that $h\circ f_p \simeq 0$, then there is a canonical
homotopy class of  maps $$H :  \cycd{0}{S^n}\otimes \bbz_p \to Y$$
satisfying $H \circ q \simeq h_*$.
}
\pf
Let $\Psi^{-1} : \cycd{0}{S^n}\otimes \bbz_p \to \cycd{0}{M(\bbz_p,n)}$
be a homotopy inverse of the canonical $\Psi$ defined in the previous
Corollary, and let $\overline{h}_* : \cycd{0}{M(\bbz_p,n)} \to Y$ be the
homomorphism established in Corollary A.2. 
Define $H \equdef \overline{h}_* \circ \Psi^{-1}$, and note that
$H\circ q = \overline{h}_* \circ \Psi^{-1} \circ q \simeq
\overline{h}_*\circ \imath_* = (\overline{h}\circ \imath)_* = h_*$,
where the first equivalence follows from Corollary A.3. This is the
desired $H$.
\qed
\bigskip

\noindent{\bf A.2.\ \ Canonical splittings}
\medskip

It is a general theorem of J. Moore that any topological abelian group
is homotopy equivalent to a product of Eilenberg-MacLane spaces.  
However, there are many inequivalent such splittings, and for the
results in these  papers and in [\LM$_1$] one makes a canonical
choice.  For the particular examples of cycle groups that we study,
the choice depends on the structure of $\bbp^n$ as a symmetric product
of  $\bbp^1$.  However, in many cases the canonical splitting is
determined  {\sl purely homotopy theoretically}.  This is the main 
result of this Appendix.  The existence of a theorem of this type was
first pointed out to the  authors by Eric Friedlander.

Throughout this discussion, the ring $R$ will always be either $\ 
\bbz\ $ or $\ \bbz/n,\ $ and we shall  use the specific model
$\cz_0(S^n)_o$  for the Eilenberg-MacLane space $K(\bbz, n)$.   
More generally for any  finitely generated module $N$ over $R$, we
shall take  $K(N, n) =  \cz_0(S^n)_o\otimes_{\bbz} N$.

Let us fix a finitely generated $R$-module $N$ and denote 
$K=K(N,n)$.
One has an isomorphism
$$
h_{K,R} : \pi_n(K)\otimes R = N\otimes R 
@>{\ \cong\ }>> H_n(K;R),
$$
obtained as the composition of isomorphisms
$N\otimes R = \pi_n(K)\otimes R 
    @>{h_K\otimes I}>>
H_n(K;\bbz)\otimes R 
    @>{\nu_K}>> H_n(K;R),$
where $h_K$ is the Hurewicz map and $\nu_K$ is provided by the 
universal
coefficients theorem for homology.

It follows from the universal coefficients theorem for cohomology 
that
$$
\Psi_{K,R}\ : \  H^n(K;N) \to \text{Hom}_R(H_n(K;R),N)
$$ is an isomorphism, and the fundamental class $\iota_n\in 
H^n(K;N)$ is
defined so that $\Psi_{K,R}(\iota_n )$ is the composition
$ H_n(K;R) @>{h^{-1}_{K,R}}>> \pi_n(K)\otimes R
= N\otimes R @>{\mu_N}>> N,$
where the latter map gives the $R$-module structure on $N$. 
Therefore,
$$
\Psi_{K,R}(\iota_n) = \mu_N\circ h^{-1}_{K,R}. \tag{A.1}
$$

We now examine these maps under the Dold-Thom theorem, which 
gives
natural isomorphisms
$$
d_{Y,R}\  : \  \pi_n(\cz_0(Y)\otimes R) \to  H_n(Y;R) \tag{A.2}
$$
for any  $CW$-complex $Y$ and for all $n$.

Since $K$ is a topological abelian group, one has a topological 
homomorphism
$$
t_K \ :\  \cz_0(K) \to K \tag{A.3}
$$
such that the composition
$
\cz_0(K)\otimes R @>{t_K\otimes I}>>
K\otimes R
@>{\mu_K}>> K
$
induces a left inverse
$$
t_{K,R} \ :\  \cz_0(K)\otimes R \to K \tag{A.4}
$$
to the natural inclusion
$$
j_{K,R} \ :\  K \to \cz_0(K)\otimes R . \tag{A.5}
$$
In the level of homotopy groups, this map fits into a commutative 
diagram
$$
\CD
N\otimes R= \pi_n(K)\otimes R @>{\cong}>{h_{K,R}}> H_n(K;R) \\
@V{\mu_N}VV  @A{\cong}A{d_{K,R}}A \\
N= \pi_n(K) @<<{\pi_n(t_{K,R})}<
\pi_n(\cz_0(K)\otimes R),
\endCD
$$
which together with  (A.1) implies that
$$
\Psi_k(\iota_n) = \mu_n\circ h^{-1}_{K,R} = \pi_n(t_{K,R}) \circ
d_{K,R}^{-1}. \tag{A.6}
$$

We now consider a map $f \ :\  Y\to K=K(N,n)$, representing a class
$[f] \in H^n(Y,N)$. From the commutative diagram
$$
\CD
H^n(Y;N) @>{\Psi_Y}>> \text{Hom}_R ( H_n(Y;R),N) @>>> 0 \\
@A{H^n(f)}AA @AA{H_n(f)^*}A @. \\
H^n(K;N) @>{\Psi_K}>> \text{Hom}_R ( H_n(K;R),N) @>>> 0 \\
\endCD
$$
one concludes that 
$$
\Psi_Y([f]) = \Psi_Y(H^n(f)(\iota_n)) = \Psi_K(\iota_n)\circ H_n(f) 
\tag{A.7}
$$
and (A.6) gives
$$
\Psi_Y([f]) = \pi_n(t_{K,R})\circ d_{K,R}^{-1} \circ H_n(f)
.
$$

Let
$$
\ov{f} \ :\ \cz_0(Y)\otimes R \to K \tag{A.8}
$$
be the $R$-module homomorphism given by the composition
$
\cz_0(Y) \otimes R 
@>{\cz_0(f)\otimes I}>>
\cz_0(K)\otimes R @>{t_{K,R}}>> K.
$
This map induces a commutative diagram
$$
\CD
\pi_n(\cz_0(Y)\otimes R) @>{\pi_n(\cz_0(f)\otimes I)}>> 
\pi_n(\cz_0(K)\otimes R) @>{\pi_n(t_{K,R})}>> \pi_n(K) = N \\
@V{d_{Y,K}}VV @V{d_{K,R}}VV @| \\
H_n(K;R) @>>{H_n(f)}> H_n(K;R) @>>{\Psi_K(\iota_n)}> \pi_n(K)=N,
\endCD
$$
where the left square commutes by the naturality of the Dold-Thom
isomorphism, and the right square commutes by (A.6).

It follows that
$
\pi_n(\ov{f}) = \Psi_K(\iota_n)\circ H_n(f)\circ d_{Y,R},
$
and by (A.7) one concludes that
$$
\pi_n(\ov{f})  = \Psi_Y([f])\circ d_{Y,R}. \tag{A.9}
$$

We now use the formulae above to prove the following result.

\Theorem {A.5}   {\sl Let $Y$ be a connected finite complex and $R= 
\bbz$ 
or $\bbz/n$. Suppose that
$$
\Psi_Y \ :\ H^n(Y;\, H_n(Y;\, R)) \ \to \ \Hom\left( H_n(Y;\, R),\,  
H_n(Y;\,R)\right)
$$
is an isomorphism for all $n$.  Then there exists a homotopy 
equivalence
$$
\cz_0(Y) \otimes_{\bbz} R \ \underset{\cong}\to {\ar {\a}} \ 
\prod_{k\geq 0} K( H_n(Y;\, R),\, n), 
$$
{\bf unique up to homotopy} with the property that:

\indent (i) \ $\a$ is a $R$-module homomorphism

\indent (ii) The composition
$$
Y \ \subset \ \cz_0(Y) \otimes R \ \ar {\a} 
\prod_{k\geq 0} K( H_n(Y;\, R),\, n)  
$$
\indent\ \ \ \ \ \ classifies the identity element in 
$H^*(Y;\, H_*(Y;\, R)) \ \cong \ \Hom\left( H_*(Y;\, R),\,  
H_*(Y;\,R)\right)$.
}

\pf
Unless otherwise indicated all homology groups have
coefficients in $R$.

We first prove existence.  For each $n$ there exists a map 
$$
f_n \ :\ Y \arr K(H_nY, n)
$$
which classifies the identity element Id $\in \Hom(H_nY, H_nY),\ $ 
i.e., $\Psi_{Y} ( [f_n] ) = \text{Id}$.
Now, define
$\ F \ :\ \cz_0(Y)\otimes R \to \prod_{n\geq 0} K(H_n(Y),n )$
by $F = \prod_{n\geq 0} \ov{f}_n$, where $\ov{f}_n$ is defined as in 
(A.8),
and note that $F$ satisfies the two conditions of the theorem. 

In the level of homotopy groups one has
$$
\pi_n(F) = \pi_n(\ov{f}_n) : \pi_n(\cz_0(Y)) \to 
\pi_n(K(H_n(Y;R),n))=H_n(Y;R),
$$
and formula (A.9) shows that
$\pi_n(f) = \Psi_Y([f_n])\circ d_{Y,R} = \text{Id}\circ 
d_{Y,R}=d_{Y,R}$, the Dold-Thom
isomorphism itself. It follows that $F$ is a homotopy equivalence.

For uniqueness suppose we are given 
$$
Y \ \underset{j_{Y,R}}\to{\hookrightarrow} \ \cz_0(Y)\otimes G\ \ 
 \underset{\underset{\b}\to {\arr}} \to {\ar{\a}} \ \ 
 \prod_n K(H_nY, n)
$$
where $\a$ and $\b$ are homotopy equivalences with properties 
(i) and (ii) above.  These properties imply that $\a\circ j$ is 
homotopic to
$\b\circ j$.  Since $\a$ and $\b$ are $R$-module homomorphsims 
and since $j_{Y,R}$ generates $\cz_0(Y)\otimes R$ as an $R$-module, 
this implies that 
$\a = \ov{\a\circ j_Y}$ is homotopic to $\b= \ov{\b\circ j_Y}$.  \qed

\bigskip

\noindent{\bf A.3.\ \ Explicit splittings of spaces of algebraic cycles}
\medskip

We have exhibited in Lemma 9.9 an explicit representative
$$\bb n \ : \ S^{2n} \to \calz^n(\bbp_{\bbc}^{2n-1})$$ for the generator
of $\pi_{2n} \cz^{\infty} \cong \bbz$.

Now, denote by $\phi_n  \ : \ S^{2n} \to  \cz^{\infty} $ the composition
$ S^{2n} \xrightarrow{\bb n} \calz^n(\bbp_{\bbc}^{2n-1})
\hookrightarrow \cz^{\infty}$, and by $\phi_0 : S^0 = \{ 0, 1 \} \to
\cz^{\infty}$ the map sending $0$ to $0$ and $1$ to $1$.  

Using the universal properties of free abelian groups on spaces, the
maps $\phi_n$ induce morphisms of abelian topological groups
$$
\Phi_n \ : \ K(\bbz,2n) = \cz_0(S^{2n})_o \to \cz^\infty.
$$
These maps can be added to give a  morphism of topological groups
$$
\Phi \ : \ \prod_{n\geq 0} K(\bbz,2n) \to \cz^\infty
$$
from the weak product $ \prod_{n\geq 0} K(\bbz,2n)$ to $ \cz^{\infty}$,
defined by $\Phi( \prod_n \ x_n ) = \sum_{n\geq 0}\ \Phi_n(x_n).$
It is immediate from its construction that this map is a homotopy
equivalence. 
\bigskip

\heading References \endheading

\bigskip

\item{[\A]} M. F. Atiyah, {\it K-theory and Reality}, Quart. J. Math. 
Oxford
(2),{\bf 17} (1966), 367-386.
   
\medskip

\item{[\BLLMM]}
Boyer, C. P., H.B. Lawson, Jr, P. Lima-Filho, B. Mann, and M.-L. 
Michelson,
{\it Algebraic cycles and infinite loop spaces,}
Inventiones   Math. , {\bf 113} (1993), 373-388.
 
\medskip

\item{[\CW]}
Costenoble, S. R.  and  S. Waner,
{\it Fixed Set Systems of Equivariant Infinite Loop Spaces,}
Trans. Amer. Math. Soc.  {\bf 326} no. 2 (1991), 485-505.

\medskip

\item{[\DT]}
Dold, A.  and R. Thom,
{\it Quasifaserungen und unendliche symmetrische produkte,}
Ann. of Math. (2) {\bf 67} (1956), 230-281.

\medskip

\item{[\dS]}
Dos Santos, P.
{\it Algebraic cycles on real varieties and $\bbz_2$-equivariant 
homotopy theory}, Ph.D. Thesis, Stony Brook, 1999.

\medskip

\item{[\D]}
Dugger, D.
{\it A Postnikov tower for algebraic K-theory,} Ph.D. Thesis, 
M.I.T., 1999.

\medskip

\item{[\FM]}
Friedlander, E. and B. Mazur,  {\it Filtrations on the homology of
algebraic varieties}, Memoire of the A. M. S., no 529,  1994.

\medskip

\item{[\FL$_1$]}
Friedlander, E.  and H.B. Lawson, Jr., {\it A theory of algebraic
cocycles}, Annals of Math., {\bf 136} (1992), 361-428.

\medskip

\item{[\FL$_2$]}
Friedlander, E.  and H.B. Lawson, Jr., {\it Duality relating spaces of
algebraic cocycles and cycles}, Topology {\bf 36} no.2 (1997),
533-565.

\medskip

\item{[\Fu]}
 Fulton, W., {\it Intersection theory},  Springer-Verlag, New York,
1984.

\medskip

\item{[\Lam]}
Lam T.-K.,
{\it Spaces of Real Algebraic Cycles and Homotopy Theory}, Ph.D. 
thesis, SUNY,
Stony Brook, 1990.

\item{[\L$_1$]}
Lawson, H.B. Jr,
{\it Algebraic cycles and homotopy theory}, Ann. of Math. {\bf 129} 
(1989),
253-291.

\medskip

\item{[\L$_2$]}
Lawson, H.B. Jr,
{\it Spaces of algebraic cycles}, pp 137-213 in ``Surveys in  
Differential
Geometry'',  vol. 2, 1995, International Press, Boston.

\medskip

\item{[\LLM$_1$]}
Lawson, H.B. Jr, P.C. Lima-Filho and M.-L. Michelsohn,
{\it Algebraic cycles and equivariant cohomology 
theories},   Proc.   London Math. Soc. (3) {\bf 73} (1996), 679-720.

\item{[\LLM$_2$]}
Lawson, H.B. Jr, P.C. Lima-Filho and M.-L. Michelsohn,
{\it On equivariant algebraic suspension}, J. Algebraic Geom. (to 
appear).

\medskip

\item{[\LLM$_3$]}
Lawson, H.B. Jr, P.C. Lima-Filho and M.-L. Michelsohn,
{\it Algebraic cycles and the classical groups, Part II}, Stony Brook 
Preprint,
1998.

\medskip

\item{[\LM$_1$]}
Lawson, H.B. Jr. and M.-L. Michelsohn,
{\it Algebraic cycles, Bott periodicity, and the Chern characteristic 
map}, in
{\it The Mathematical Heritage of Hermann Weyl, } A.M.S., 
Providence, 1988, pp.
241-264.

\medskip

\item{[\LM$_2$]}
Lawson, H.B. Jr and M.-L. Michelsohn,
{\it Algebraic cycles and group actions} in {\it Differential Geometry}, 
Longman
Press, 1991, pp. 261-278.

\medskip

\item{[\LMS]}
Lewis, Jr., L. G. and J. P. May and M. Steinberger,
{\it Equivariant Stable Homotopy Theory}, LNM, vol. 1213,
  Springer-Verlag, New York, NY, 1986.

\item{[\Li$_1$]}
 Lima-Filho, P.C.
{\it Lawson homology for quasiprojective varieties,} 
Compositio Math {\bf 84} (1992), 1-23.
\medskip

\item{[\Li$_2$]}
 Lima-Filho, P.C.,
{\it Completions and fibrations for topological monoids,}
Trans. Amer. Math. Soc.  {\bf 340} no. 1 (1993), 127-146.
\medskip

\item{[\Li$_3$]}
Lima-Filho, P.C.,  
{\it On the generalized cycle map,}
J. Diff. Geom. {\bf 38} (1993), 105-130.

\medskip

\item{[\Li$_4$]}
 Lima-Filho, P.C.,
{\it On the equivariant homotopy of free abelian groups on $G$-
spaces and $G$-spectra,}
Math. Z. {\bf 224} (1997), 567-601.

\medskip

\item{[\Li$_5$]}
 Lima-Filho, P.C.,
{\it On the topological group structure of algebraic cycles,}
Duke Math. J. {\bf 75} (1994), no. 2,467--491.

\medskip

\item{[\May$_1$]}
 May, J.~P. {\it The geometry of iterated loop spaces}, LNM, vol. 271,
  Springer-Verlag, New York, NY, 1972.

\medskip

\item{[\May$_2$]}
May, J. P.,
{\it Classifying Spaces and Fibrations,}
Mem. Amer. Math. Soc. {\bf 155}, 1975.

\medskip

\item{[\May$_3$]}
May, J. P., {\it ${E}_\infty$-ring spaces and ${E}_\infty$-ring 
spectra}, LNM,
  vol. 577, Springer-Verlag, New York, 1977.

\medskip

\item{[\May$_4$]}
May, J. P., et al. {\it Equivariant homotopy and cohomology
theory}, CBMS vol. 91, American Mathematical  Society, 1996.

\medskip

\item{[\MS]}
Milnor, J. and Stasheff, J. D., {\it Characteristic Classes}, Annals of 
Math.
Studies no. 76, Princeton Univ. Press, 1974.

\medskip

\item{[\Mo$_1$]}
Mostovoy, J.,
{\it Spaces of real algebraic cycles on $\bbp^n$,}
Russian Math. Surveys {\bf 53} (1998), no.1 (319)

\medskip

\item{[\Mo$_2$]}
Mostovoy, J.,
{\it Algebraic cycles and anti-holomorphic involutions on projective spaces,}
UNAM Preprint, 1998. 

\medskip

\item{[\Se]}
Segal, G. {\it The multiplicative group of classical cohomology}, Quart. 
J.
  Math. Oxford Ser. (2) {\bf 26} (1975), 289--293.

\enddocument

\enddocument